\documentclass[12pt]{amsart}
\usepackage{amsmath,amssymb,amsfonts}

\usepackage{epsf}

\usepackage{epic}
\usepackage{eepic}

\def\tagref#1{(#1)}
\def\proclaim#1{\vskip12pt\bgroup\noindent{\sc#1}.\enskip\it\ignorespaces}
\def\endproclaim{\egroup\vskip6pt}
\renewcommand{\tfrac}[2]{{\textstyle \frac{#1}{#2}}}

\DeclareSymbolFont{AMSb}{U}{msb}{m}{n}
\DeclareMathSymbol{\C}{\mathalpha}{AMSb}{"43}
\DeclareMathSymbol{\F}{\mathalpha}{AMSb}{"46}
\DeclareMathSymbol{\N}{\mathalpha}{AMSb}{"4E}
\DeclareMathSymbol{\Q}{\mathalpha}{AMSb}{"51}
\DeclareMathSymbol{\R}{\mathalpha}{AMSb}{"52}
\DeclareMathSymbol{\Z}{\mathalpha}{AMSb}{"5A}
\DeclareMathSymbol{\nmid}{\mathrel}{AMSb}{"2D}
\def\nutwid{\overset {\text{\lower 3pt\hbox{$\sim$}}}\nu}
\def\qtwid{\overset {\text{\lower 3pt\hbox{$\sim$}}}q}
\def\pitwid{\overset {\text{\lower 3pt\hbox{$\sim$}}}\pi}
\def\Etwid{\overset {\text{\lower 3pt\hbox{$\sim$}}}E}
\def\etwid{\overset {\text{\lower 3pt\hbox{$\sim$}}}e}
\def\pio{\overline{\pi}}
\def\mycap#1#2{\hbox to \textwidth{\hfil{\textsc{Fig #1}.{\ #2}}\hfil}}
\def\qbin#1#2{\left[
\matrix #1 \\
#2
\endmatrix
\right]_q}
\def\qbinb#1#2{\left[
\matrix #1 \\
#2
\endmatrix
\right]_{q^2}}
\def\rank{\mbox{rank}}
\def\crank{\mbox{crank}}
\def\sign{\mbox{sign}}
\def\Chat{\widehat{C}}
\def\Ghat{\widehat{G}}
\def\chat{\widehat{c}}
\def\Ehat{\widehat{E}}
\def\CT{\mbox{CT}}
\def\2phi1#1#2#3{
{}_{2}\phi_1\left(\begin{array}{l} #1\\ 
                                   #2 \end{array}; #3\right)}
\def\qinf#1{(#1)_\infty}
\def\SPC{\mbox{SPC}}
\def\au{\mbox{au}\,}
\def\LHS{\mbox{LHS}}
\begin{document}

\renewcommand{\theequation}{\arabic{section}.\arabic{equation}}

\title{Some Observations on Dyson's\\
New Symmetries of Partitions}

\author{Alexander Berkovich}
\email{alexb\char"40math.ufl.edu} 
\address{Department of Mathematics, University of Florida,
Gainesville, FL 32611}

\author{Frank G. Garvan}
\email{frank\char"40math.ufl.edu} 
\address{Department of Mathematics, University of Florida,
Gainesville, FL 32611}

\thanks{Research of the second author supported in part by the NSF under grant 
number DMS-9870052.}
\date{April 17, 2002}
%%\date{April 11, 2002}
%%\date{April 9, 2002}
%%\date{April 8, 2002}
%%\date{March 24, 2002}
%%\date{March 12, 2002}

\begin{abstract}
We utilize Dyson's concept of the adjoint of a partition to
derive an infinite family of new polynomial analogues of Euler's Pentagonal 
Number Theorem.
We streamline Dyson's bijection relating partitions with crank $\le k$
and those with $k$ in the Rank-Set of partitions. Also, we extend
Dyson's adjoint of a partition to MacMahon's ``modular'' partitions
with modulus $2$. This way we find a new combinatorial proof of
Gauss's famous identity.
We give a direct combinatorial proof that for $n>1$ the partitions of $n$
with crank $k$ are equinumerous with partitions of $n$ with
crank $-k$.
\end{abstract}

\maketitle
\pagestyle{myheadings}
\markboth{A. Berkovich and F.G. Garvan}{Dyson's symmetries of partitions}

\setcounter{equation}{0}
\section{Introduction}
\label{sec:intro}

%%%%%%SECTION 1
Let $p(n)$ denote the number of unrestricted partitions of $n$.
Ramanujan discovered three beautiful arithmetic properties of $p(n)$,
namely:
\begin{align}
p(5n+4) &\equiv 0 \pmod{5},\tag{1.1}\\
p(7n+5) &\equiv 0 \pmod{7},\tag{1.2}\\
p(11n+6) &\equiv 0 \pmod{11}.\tag{1.3}
\end{align}
The partition congruences modulo $5$ and $7$ were proved by Ramanujan 
in \cite{Ra1}. In \cite{Ra2} he proved \tagref{1.3} by a different method.
The most elementary proof of \tagref{1.3} similar to the one in \cite{Ra1}
is due to Winquist \cite{Win}.

Dyson \cite{D1} discovered empirically remarkable combinatorial
interpretations of \tagref{1.1} and \tagref{1.2}.
Defining the {\it rank} of a partition as the {\it largest part minus
the number of parts}, he observed that
\begin{align}
N(k,5,5n+4) &= \frac{p(5n+4)}{5},\quad 0 \le k \le 4, \tag{1.4}\\
N(k,7,7n+5) &= \frac{p(7n+5)}{7},\quad 0 \le k \le 6, \tag{1.5}
\end{align}
where $N(k,m,n)$ denotes the number of partitions of $n$ with
rank congruent to $k$ modulo $m$.
Identities \tagref{1.4} and \tagref{1.5} were later proved by
Atkin and Swinnerton-Dyer \cite{ASD}.
However, the rank failed to explain \tagref{1.3}, and so Dyson conjectured
the existence of some analogue of the rank that would explicate the
Ramanujan congruence modulo $11$. He named his hypothetical statistic
the {\it crank}.

Forty four years later, Andrews and Garvan \cite{AG}, building on the
work of Garvan \cite{Ga1}, finally unveiled
Dyson's crank of a partition $\pi$:
\begin{equation}
\mbox{crank}(\pi) =
\begin{cases}
  \lambda(\pi), &\mbox{if $\mu(\pi)=0$}, \\
  \nutwid(\pi) - \mu(\pi), &\mbox{if $\mu(\pi)>0$},
\end{cases}
\tag{1.6}
\end{equation}
where $\lambda(\pi)$ denotes the largest part of $\pi$,
$\mu(\pi)$ denotes the number of ones in $\pi$ and $\nutwid(\pi)$
denotes the number of parts of $\pi$ larger than $\mu(\pi)$.

Remarkably, the crank provides combinatorial interpretations
of all three Ramanujan congruences \tagref{1.1}--\tagref{1.3}.
Namely,
\begin{align}
M(k,5,5n+4) &= \frac{p(5n+4)}{5},\quad 0 \le k \le 4, \tag{1.7}\\
M(k,7,7n+5) &= \frac{p(7n+5)}{7},\quad 0 \le k \le 6, \tag{1.8}\\
M(k,11,11n+6) &= \frac{p(11n+6)}{11},\quad 0 \le k \le 10, \tag{1.9}  
\end{align}
where $M(k,m,n)$ denotes the number of partitions of $n$ with
crank congruent to $k$ modulo $m$.

Let $P_m(q)$ denote the generating function
\begin{equation}
P_m(q) = \sum_{n=1}^\infty p_m(n) q^n,
\tag{1.10}
\end{equation}
where $p_m(n)$ is the number of partitions of $n$ with rank $m$.
Here we are using the convention that $p_m(0)=0$.
As a practical tool for his empirical calculations Dyson used the 
following formula for $P_m(q)$:
\begin{equation}
P_m(q) = \frac{1}{(q)_\infty}
\sum_{j\ge1} (-1)^{j-1} (1-q^j) q^{\frac{j(3j-1)}{2}+ |m|j},
\tag{1.11}
\end{equation}
with
\begin{equation}
(q)_\infty = \prod_{j\ge1}(1-q^j),\quad \mbox{for $|q|<1$.}
\tag{1.12}
\end{equation}
For later use we also define
\begin{equation}
(a;q)_n = (a)_n = \prod_{j=0}^{n-1} (1 - a q^j), \quad n\ge 0
\tag{1.13}
\end{equation}
and note that $\frac{1}{(q)_\infty}$ is the generating function
for unrestricted partitions.

Dyson knew how to prove \tagref{1.11} in 1942 \cite{D2}.
However the first published proof of \tagref{1.11} was given
by Atkin and Swinnerton-Dyer \cite{ASD} in 1954.
In 1968, Dyson \cite{D3} found a simple combinatorial argument
which not only explained \tagref{1.11} but also led to a new proof of Euler's
celebrated pentagonal number theorem:
\begin{equation}
1 = \frac{1}{(q)_\infty} \sum_{j=-\infty}^{\infty} (-1)^j q^{\frac{j(3j-1)}{2}}.
\tag{1.14}
\end{equation}

To paraphrase Dyson's argument in  \cite{D3} we introduce the
generating function
\begin{equation}
Q_m(q) = \sum_{n\ge1} \qtwid_m(n) q^n,
\tag{1.15}
\end{equation}
where $\qtwid_m(n)$ is the number of partitions of $n$
with rank $\ge m$. We adopt the convention that
\begin{equation}
\qtwid_m(0)=0.
\tag{1.16}
\end{equation}
Clearly,
\begin{equation}
P_m(q) = Q_m(q) - Q_{m+1}(q).
\tag{1.17}
\end{equation}
Next, following the treatment in \cite{D3} we will show that
\begin{equation}
Q_m(q) + Q_{1-m}(q) + 1 = \frac{1}{(q)_\infty},
\tag{1.18}
\end{equation}
and for $m\ge0$
\begin{equation}
Q_m(q) = q^{m+1}\left(Q_{-2-m}(q) + 1\right).
\tag{1.19}
\end{equation}

To prove \tagref{1.18} we note that any given nonempty partition $\pi$
counted by $\frac{1}{(q)_\infty}$ has either rank $\ge m$ or
rank $< m$. If rank $\ge m$, then $\pi$ is counted by $Q_m(q)$.
If rank $< m$, then we conjugate $\pi$ to get $\pi^{*}$ as
illustrated in Fig.\ 1.

\begin{figure}[h!]
\begin{center}
\setlength{\unitlength}{0.00033333in}
\renewcommand{\dashlinestretch}{30}
\begin{picture}(15473,9927)(0,-10)
%%\begin{picture}(17473,9927)(0,-10)
%%\begin{picture}(40473,9927)(0,-10)
%%\begin{picture}(5000,40000)(0,-10)
\put(5392.500,5080.500){\arc{597.369}{5.6045}{8.3513}}
\blacken\path(5373.560,4812.271)(5250.000,4818.000)(5356.328,4754.799)(5330.461,4793.875)(5373.560,4812.271)
\path(1125,9318)(7125,9318)
\path(1125,9318)(1125,4518)
\path(2325,9393)(2325,9243)
\path(3600,9393)(3600,9243)
\path(4725,9393)(4725,9243)
\path(5925,9393)(5925,9243)
\path(1050,8118)(1200,8118)
\path(1050,6918)(1200,6918)
\path(1050,5718)(1200,5718)
\path(1050,8718)(1200,8718)
\path(1050,7518)(1200,7518)
\path(1050,6318)(1200,6318)
\path(1050,5118)(1200,5118)
\path(3000,9393)(3000,9243)
\path(1725,9393)(1725,9243)
\path(4125,9393)(4125,9243)
\path(5325,9393)(5325,9243)
\path(6525,9393)(6525,9243)
\dashline{60.000}(1125,9318)(5925,4518)
\drawline(2475,4443)(2475,4443)
\path(975,9318)(900,9168)(900,6993)
	(825,6918)(900,6843)(900,4668)(975,4518)
\path(1125,9468)(1275,9543)(1275,9543)
	(4050,9543)(4125,9618)(4200,9543)
	(6975,9543)(7125,9468)
\path(1125,4518)(1725,4518)(1725,5718)
	(3525,5718)(3525,6318)(4125,6318)
	(4125,6918)(4725,6918)(4725,7518)
	(5925,7518)(5925,8718)(7125,8718)(7125,9318)
\path(11464,9322)(11464,3323)
\path(11464,9322)(16262,9322)
\path(11388,8122)(11538,8122)
\path(11389,6847)(11539,6847)
\path(11389,5723)(11539,5723)
\path(11389,4522)(11539,4522)
\path(12663,9397)(12663,9247)
\path(13864,9397)(13864,9247)
\path(15063,9397)(15062,9247)
\path(12065,9397)(12065,9247)
\path(13265,9399)(13265,9248)
\path(14465,9397)(14465,9247)
\path(15665,9398)(15665,9248)
\path(11388,7448)(11539,7447)
\path(11389,8721)(11539,8721)
\path(11389,6321)(11539,6321)
\path(11388,5121)(11539,5122)
\path(11388,3922)(11540,3921)
\dashline{60.000}(11464,9322)(16264,4521)
\drawline(16337,7973)(16337,7973)
\path(11464,9472)(11614,9547)(13789,9547)
	(13864,9622)(13939,9547)(16114,9547)(16264,9472)
\path(11314,9322)(11238,9172)(11238,9172)
	(11239,6396)(11163,6321)(11238,6246)
	(11238,3472)(11314,3323)
\path(16262,9322)(16264,8721)(15064,8722)
	(15062,6922)(14465,6922)(14465,6321)
	(13864,6321)(13864,5723)(13265,5722)
	(13265,4523)(12065,4522)(12064,3321)(11464,3323)
\thicklines
\path(6675,6468)(8325,6468)
\blacken\thinlines
\path(8205.000,6378.000)(8325.000,6468.000)(8205.000,6558.000)(8241.000,6468.000)(8205.000,6378.000)
\put(-1000,6993){$\pi:$}
\put(-100,6343){$\nu(\pi)$}
\put(3900,9768){$\lambda(\pi)$}
\put(10200,6243){$\lambda(\pi)$}
\put(13725,9768){$\nu(\pi)$}
\put(9350,6993){$\pi^{*}:$}
\end{picture}
\end{center}
%%\caption{{Fig. 1}\enskip Conjugation of a partition $\pi$}
\vskip -60pt
%%\mycap{1}{\enskip Conjugation of a partition $\pi$}
\textsc{Fig 1}.\ \parbox[t]{5in}{
Conjugation of a partition $\pi$ with largest part $\lambda(\pi)$
and the number of parts $\nu(\pi)$, $\lambda(\pi)-\nu(\pi)<m$.}
\end{figure}

It is obvious in this case
that $\mbox{rank}(\pi^{*})\ge1-m$. Hence, $\pi^{*}$ is
counted by $Q_{1-m}(q)$. Finally, the empty partition is counted by $1$
on the left side of \tagref{1.18} and on the right side
by $\frac{1}{(q)_\infty}$.
The proof of \tagref{1.19} is more subtle. Here we will use a different
conjugation transformation (Dyson's adjoint) as follows.
Consider some partition $\pi$ with $\mbox{rank}(\pi)\ge m\ge0$.
This partition is counted by $Q_m(q)$ in \tagref{1.19}.
Clearly,
\begin{equation}
\lambda(\pi) - \nu(\pi) \ge m \ge 0.
\tag{1.20}
\end{equation}
Let us now remove the largest part of $\pi$ to end up with $\pitwid$.
Next, we conjugate $\pitwid$ to get $\pitwid^{*}$.
Finally, we attach to $\pitwid^{*}$ a new largest part of
size $\lambda(\pi) - m - 1$.
These transformations are illustrated in Fig.\ 2.

\begin{figure}[h!]
\begin{center}
\setlength{\unitlength}{0.00049167in}
\begingroup\makeatletter\ifx\SetFigFont\undefined%
\gdef\SetFigFont#1#2#3#4#5{%
  \reset@font\fontsize{#1}{#2pt}%
  \fontfamily{#3}\fontseries{#4}\fontshape{#5}%
  \selectfont}%
\fi\endgroup%
{\renewcommand{\dashlinestretch}{30}
\begin{picture}(8112,1239)(0,-10)
\path(5400,987)(5700,987)(5700,687)
	(5400,687)(5400,987)
\path(5400,687)(5700,687)(5700,387)
	(5400,387)(5400,687)
\path(5400,387)(5700,387)(5700,87)
	(5400,87)(5400,387)
\path(5700,987)(6000,987)(6000,687)
	(5700,687)(5700,987)
\path(450,912)(750,912)(750,612)
	(450,612)(450,912)
\path(750,912)(1050,912)(1050,612)
	(750,612)(750,912)
\path(1050,912)(1350,912)(1350,612)
	(1050,612)(1050,912)
\path(1350,912)(1650,912)(1650,612)
	(1350,612)(1350,912)
\thicklines
\path(1650,312)(2250,312)
\blacken\thinlines
\path(2130.000,222.000)(2250.000,312.000)(2130.000,402.000)(2166.000,312.000)(2130.000,222.000)
\path(450,612)(750,612)(750,312)
	(450,312)(450,612)
\path(750,612)(1050,612)(1050,312)
	(750,312)(750,612)
\path(1050,612)(1350,612)(1350,312)
	(1050,312)(1050,612)
\path(450,312)(750,312)(750,12)
	(450,12)(450,312)
\path(3000,687)(3300,687)(3300,387)
	(3000,387)(3000,687)
\path(3300,687)(3600,687)(3600,387)
	(3300,387)(3300,687)
\path(3600,687)(3900,687)(3900,387)
	(3600,387)(3600,687)
\path(3000,387)(3300,387)(3300,87)
	(3000,87)(3000,387)
\thicklines
\path(4050,387)(4650,387)
\blacken\thinlines
\path(4530.000,297.000)(4650.000,387.000)(4530.000,477.000)(4566.000,387.000)(4530.000,297.000)
\thicklines
\path(6150,387)(6750,387)
\blacken\thinlines
\path(6630.000,297.000)(6750.000,387.000)(6630.000,477.000)(6666.000,387.000)(6630.000,297.000)
\path(7500,912)(7800,912)(7800,612)
	(7500,612)(7500,912)
\path(7500,612)(7800,612)(7800,312)
	(7500,312)(7500,612)
\path(7500,312)(7800,312)(7800,12)
	(7500,12)(7500,312)
\path(7800,912)(8100,912)(8100,612)
	(7800,612)(7800,912)
\path(7500,1212)(7800,1212)(7800,912)
	(7500,912)(7500,1212)
\path(7800,1212)(8100,1212)(8100,912)
	(7800,912)(7800,1212)
\put(0,360){$\pi$:}
\put(2550,360){$\pitwid$:}
\put(4950,360){$\pitwid^{*}$:}
\put(7050,360){$\pi'$:}
%\put(0,387){$\pi$:}
%\put(2550,312){$\piwtid$:}
%\put(4950,312){$\pitwid^{*}$:}
%\put(7050,387){$\pi'$:}
\end{picture}
}
\end{center}
\mycap{2}{\enskip Dyson's adjoint of $\pi=4+3+1$ with $\mbox{rank}=m=1$.}
\end{figure}

Note that the map $\pi \to \pi'$ is reversible. It is obvious that
\begin{align}
\lambda(\pi') &= \lambda(\pi) - m - 1,\tag{1.21}\\
\nu(\pi') &\le \lambda(\pi)  + 1,\tag{1.22}\\
|\pi'| &= |\pi|  -m - 1,\tag{1.23}
\end{align}
where $|\pi|$ denotes the sum of parts of $\pi$.

Since $\mbox{rank}(\pi')\ge -2-m$, we see that $\pi'$ is counted
by $Q_{-2-m}(q)$ in \tagref{1.19}, provided $|\pi'|\ne0$.
If $|\pi|=m+1$ and $\mbox{rank}(\pi)\ge m$, then
$\nu(\pi)=1$ and $\lambda(\pi)=1+m$. So in this case $\pi'$ represents
the empty partition, which is counted by $1$ in \tagref{1.19}. This
concludes the proof of \tagref{1.19}.

Combining \tagref{1.18} and \tagref{1.19}, we see that for $m\ge 0$
\begin{equation}
Q_m(q) + q^{m+1} Q_{m+3} = \frac{q^{m+1}}{(q)_\infty}.
\tag{1.24}
\end{equation}
Iterating \tagref{1.24} we obtain for $m\ge0$
\begin{align}
Q_m(q) &= \frac{q^{m+1}}{(q)_\infty} - q^{m+1} Q_{m+3} \nonumber\\
&= \frac{q^{m+1} - q^{2m+5}}{(q)_\infty} + q^{2m+5} Q_{m+6}= \cdots
\tag{1.25}\\
& = 
 \frac{1}{(q)_\infty} \sum_{j\ge1} (-1)^{j-1} q^{\frac{j(3j-1)}{2} + mj}.
\nonumber
\end{align}

We observe that \tagref{1.25} together with \tagref{1.18} with $m=0$
yields Euler's theorem \tagref{1.14}. On the other hand, \tagref{1.25}
together with \tagref{1.17} yields \tagref{1.11} with $m\ge0$.
To treat the $m<0$ case in \tagref{1.17} we make use of
\begin{equation}
P_m(q) = P_{-m}(q),
\tag{1.26}
\end{equation}
which is a straightforward consequence of the conjugation transformation.

In \cite{A1} Andrews utilized Dyson's adjoint to give
a new proof of a partition theorem due to Fine.
This seems to be the only known application of the Dyson transformation.

In the next section we will show that Dyson's formulas \tagref{1.18}, 
\tagref{1.19} can be generalized to yield a binary tree of
polynomial analogues of \tagref{1.14}. This tree contains 
Schur's well-known formula
\begin{equation}
1=\sum_{j=-\infty}^{\infty} (-1)^j q^{\frac{j(3j-1)}{2}}
\left[
\matrix 2L \\
L + \lfloor \tfrac{3j}{2} \rfloor
\endmatrix
\right]_q,
\tag{1.27}
\end{equation}
as well as a new polynomial version of \tagref{1.14}
\begin{equation}
1=\sum_{j=-\infty}^{\infty} (-1)^j q^{\frac{j(3j+1)}{2}}
\left[
\matrix 2L - j\\
L + j
\endmatrix
\right]_q,
\tag{1.28}
\end{equation}
where $\lfloor x \rfloor$ denotes the integer part of $x$ and
$q$-binomial coefficients are defined as
\begin{equation}
\left[
\matrix n + m \\
n
\endmatrix
\right]_q
=
\begin{cases}
  \frac{(q)_{n+m}}{(q)_n (q)_m}, &\mbox{if $n$, $m\ge0$}, \\
  0, &\mbox{otherwise.}       
\end{cases}
\tag{1.29}
\end{equation}

Actually, \tagref{1.27} and \tagref{1.28} are special cases
of the following more general formula
\begin{equation}
1 - \delta_{\sigma,-1} =
\sum_{j=-\infty}^{\infty} (-1)^j q^{\frac{j(3j-1)}{2} + \sigma j}
\qbin{2L + \sigma - \au(n,j)}{L + \sigma + \lfloor \frac{n+1}{n}j\rfloor},
\tag{1.30}
\end{equation}
where
$n=1$, $2$, $3$, $4$, $5$, \dots; $\sigma=-1$, $0$, $1$, and
\begin{equation}
\au(n,j) = 
\begin{cases}
   -j, &\mbox{if $n=1$},\\
   0, &\mbox{if $n=2$},\\
   \displaystyle\sum_{k=1}^{n-2}\left\lfloor \frac{j+k}{n} \right\rfloor, &\mbox{if $n>2$.}
\end{cases}
\tag{1.31}
\end{equation}
It is easy to verify that
\begin{equation}
\lim_{n\to\infty} \au(n,j) = j-1,\quad \text{if $j>0$,}
\tag{1.32}
\end{equation}
and
\begin{equation}
\au(n,-j) = -\au(n,j),\quad \text{for $n\ge1$.}
\tag{1.33}
\end{equation}

Note that \tagref{1.27} is \tagref{1.30} with $n=2$, $\sigma=0$
and \tagref{1.28} is \tagref{1.30} with $n=1$, $\sigma=0$.

In \S\S 3 and 4 we will streamline and generalize Dyson's treatment of 
partitions
with crank $\le k$. In \S5 we will use modular representations
with modulus $2$
of partitions 
in which odd parts do not repeat,
and an appropriate modification of Dyson's
adjoint transformation, to obtain a new proof
of the Gauss formula
\begin{equation}
\frac{(q^2;q^2)_\infty}{(q;q^2)_\infty}
= \sum_{j\ge0}  q^{\frac{j(j+1)}{2}}.
\tag{1.34}
\end{equation}
%%We conclude in \S6 with a brief description of some open questions for
%%future investigation.
%%%%%%%%%%%%%%%%%%%%%%%%%%%%%%%%%%%%%%%%%%%%%%%%%%%%%%%%%%%%%%%%%%%%%%%%%
%%INSERT #1
We remark that the first combinatorial proof of \tagref{1.34}
was given by Andrews in \cite{X}. This early proof uses a
Franklin-type involution and is quite different from the one given
in \S5.
%%%%%%%%%%%%%%%%%%%%%%%%%%%%%%%%%%%%%%%%%%%%%%%%%%%%%%%%%%%%%%%%%%%%%%%%%
\S6 contains a brief description of some open questions for
future research.  In Appendix A we introduce
a new type of partition transformation, termed pseudo-conjugation,
in order to prove 
directly  that for $n>1$ the partitions of $n$
with crank $k$ are equinumerous with partitions of $n$ with
crank $-k$.
We also show that self-pseudo-conjugate partitions of $n$
(introduced there) are equinumerous with partitions of $n$ into
distinct odd parts.
%%%%%%%%%%%%%%%%%%%%%%%%%%%%%%%%%%%%%%%%%%%%%%%%%%%%%%%%%%%%%%%%%%%%%%%%%
%%INSERT #2
Finally, in Appendix B we outline an alternative proof of the formula 
\tagref{5.17}. This proof was communicated to us by George Andrews \cite{Y}.
%%%%%%%%%%%%%%%%%%%%%%%%%%%%%%%%%%%%%%%%%%%%%%%%%%%%%%%%%%%%%%%%%%%%%%%%%
%%%%%END SECTION 1

\setcounter{equation}{0}
\section{Polynomial analogues of Euler's pentagonal number theorem}
\label{sec:2}
%%%%%%SECTION 2

We say that a partition $\pi$ is in the box $[L,M]$ if its largest part
does not exceed $L$ and the number of parts does not exceed $M$. In other words,
\begin{align*}
\lambda(\pi) &\le L,\\
\nu(\pi) &\le M.
\end{align*}
It is well known \cite{A2} that the generating functions for partitions
in the box $[L,M]$
is
$$
\left[
\matrix L + M \\
L
\endmatrix
\right]_q.
$$
Let us define $Q_m^L(q)$ as
\begin{equation}
Q_m^L(q) = \sum_{n\ge1} \qtwid_m^L(n) q^n,
\tag{2.1}
\end{equation}
where $\qtwid_m^L(n)$ is the number of partitions of $n$ with
rank $\ge m$ and largest part $\le L$. As before, we assume that
$\qtwid_m^L(0)=0$. Clearly, $Q_m^L(q)=0$ whenever $L \le m$. We now prove
that
\begin{equation}
Q_m^L(q) - Q_m^{L-1}(q) = q^L \qbin{2L -m -1}{L}.
\tag{2.2}
\end{equation}
To this end we observe that the left side of \tagref{2.2} counts
partitions with rank $\ge m$ and $\lambda(\pi)=L$.
We note that these partitions are in the box $[L,L-m]$. 
If we remove the largest part $L$ from one of those partitions
we obtain a partition in the
box $[L, L-m-1]$, and this partition is counted by the $q$-binomial coefficient
on the right side of \tagref{2.2}, as desired.

We now move on to derive the bounded analogues of \tagref{1.18}, \tagref{1.19},
namely:
\begin{equation}
Q_m^L(q) + Q_{1-m}^{L-m}(q) + 1 = \qbin{2L-m}{L},\qquad L >m,
\tag{2.3}
\end{equation}
and 
\begin{equation}
Q_m^L(q) = q^{m+1}\left(Q_{-2-m}^{L-1-m}(q) + 1\right),\qquad L>m\ge0.
\tag{2.4}
\end{equation}
To prove \tagref{2.3} we note that any given nonempty partition
$\pi$, counted by $\qbin{2L-m}{m}$, has either rank $\ge m$
or rank $<m$. Moreover, $\pi$ is in the box $[L,L-m]$.
Now if $\rank(\pi)\ge m$, then $\pi$ is counted by $Q_m^L(q)$.
If $\rank(\pi)< m$, then we conjugate $\pi$ to get $\pi^{*}$.
Obviously, $\pi^{*}$ is counted by $Q_{1-m}^{L-m}(q)$. Finally,
the empty partition is counted by $1$ and by the $q$-binomial coefficient
on the left and right sides of \tagref{2.3}, respectively.

Next, to prove \tagref{2.4} we observe that any partition $\pi$ counted
by $Q_m^L(q)$ is in the box $[L,L-m]$.  Performing Dyson's
transformation $\pi \to \pi'$, as explained in the previous section,
we see that $\lambda(\pi')\le L-1-m$
and $\rank(\pi')\ge -2 -m$. Therefore, if $|\pi| \ne m+1$, then
$\pi'$ is counted by $Q_{-2-m}^{L-1-m}(q)$. If $|\pi| = m+1$,
then $\pi'$ is empty. In this case it is counted by $1$ on the right side
of \tagref{2.4}.

%
% L -> L+2, m -> m+3 in (2.3) and use (2.4)
%
Combining \tagref{2.3} and \tagref{2.4} yields
\begin{equation}
Q_m^L(q) + q^{m+1} Q_{m+3}^{L+2}(q) 
= q^{m+1} \qbin{2L-m+1}{L+2}, \quad m\ge0.
\tag{2.5}
\end{equation}
We remark that when $L\le m$ the above formula becomes $0+0=0$,
and when $L$ tends to infinity \tagref{2.5} reduces to \tagref{1.24}.
Actually, it is possible to derive another bounded analogue of \tagref{1.24}.
To this end we employ \tagref{2.2} together with the
well known recurrence
\begin{equation}
\qbin{n+m}{n} = q^n \qbin{n+m-1}{n} + \qbin{n+m-1}{m-1}
\tag{2.6}
\end{equation}
to transform \tagref{2.5} as
\begin{align}
Q_m^L(q) + q^{m+1} Q_{m+3}^{L+1}(q)
&= q^{m+1}\left\{
\qbin{2L-m+1}{L+2} - q^{L+2} \qbin{2L-m}{L+2}\right\}\tag{2.7}\\
&= q^{m+1} \qbin{2L-m}{L+1},\quad m\ge 0.
\nonumber
\end{align}

The power of \tagref{2.5} and \tagref{2.7} lies in the fact that these
transformations can be employed to generate an infinite binary tree
of representations for $Q_m^L(q)$.
First we consider four special cases, namely:
\begin{align}
Q_m^L(q) &= \sum_{j\ge1} (-1)^{j-1} q^{\frac{j(3j-1)}{2} + mj}
\qbin{2L-m+j}{L-m-j}, \quad m\ge0, \tag{2.8}\\
Q_m^L(q) &= \sum_{j\ge1} (-1)^{j-1} q^{\frac{j(3j-1)}{2} + mj}
\qbin{2L-m-j+1}{L+j}, \quad m\ge0, \tag{2.9}\\
Q_m^L(q) &= \sum_{j\ge1} (-1)^{j-1} q^{\frac{j(3j-1)}{2} + mj}
\qbin{2L-m+1}{L-\lfloor-\tfrac{3j}{2} \rfloor}, \quad m\ge0, \tag{2.10}\\
Q_m^L(q) &= \sum_{j\ge1} (-1)^{j-1} q^{\frac{j(3j-1)}{2} + mj}
\qbin{2L-m}{L+\lfloor\tfrac{3j}{2} \rfloor}, \quad m\ge0. 
\tag{2.11}
\end{align}
To derive \tagref{2.8}--\tagref{2.11} we use the
iteration schemes which we denote symbolically as
\begin{align}
&\tagref{2.5} - \tagref{2.5} - \tagref{2.5} - \tagref{2.5} - \tagref{2.5} - \tagref{2.5} - \cdots, \tag{2.12}\\
&\tagref{2.7} - \tagref{2.7} - \tagref{2.7} - \tagref{2.7} - \tagref{2.7} - \tagref{2.7} - \cdots, \tag{2.13}\\
&\tagref{2.5} - \tagref{2.7} - \tagref{2.5} - \tagref{2.7} - \tagref{2.5} - \tagref{2.7} - \cdots, \tag{2.14}\\
&\tagref{2.7} - \tagref{2.5} - \tagref{2.7} - \tagref{2.5} - \tagref{2.7} - \tagref{2.5} - \cdots, \tag{2.15}
\end{align}
respectively.
For example, the scheme \tagref{2.12} means each transformation uses
only equation \tagref{2.5}, and the scheme \tagref{2.14} means
that we use both \tagref{2.5} and \tagref{2.7} in an alternating
fashion with \tagref{2.5} being used first.

Now, \tagref{2.3} with $m=0$ yields
\begin{equation}
Q_0^L(q) + Q_{1}^{L}(q) + 1 = \qbin{2L}{L}.
\tag{2.16}
\end{equation}
Equation \tagref{1.28} then follows by using \tagref{2.8} with $m=0$
and \tagref{2.9} with $m=1$. 

Schur's formula \tagref{1.27} follows in a similar fashion.
We use \tagref{2.16}, \tagref{2.10} with $m=1$, \tagref{2.11}
with $m=0$ and the fact that 
\begin{equation}
\qbin{2L}{L+a}=\qbin{2L}{L-a}.
\tag{2.17}
\end{equation}

%%%%%INSERT II
To prove \tagref{1.30} we need to consider the following periodic
iterations:
\begin{equation}
\underbrace{
\tagref{2.7} - \tagref{2.7} - \cdots - \tagref{2.7} - \tagref{2.5}}_n
-
\underbrace{
\tagref{2.7} - \tagref{2.7} - \cdots - \tagref{2.7} - \tagref{2.5}}_n
- \cdots,
\tag{2.18}
\end{equation}
and
\begin{equation}
\underbrace{
\tagref{2.5} - \tagref{2.5} - \cdots - \tagref{2.5} - \tagref{2.7}}_n
-
\underbrace{
\tagref{2.5} - \tagref{2.5} - \cdots - \tagref{2.5} - \tagref{2.7}}_n
- \cdots.
\tag{2.19}
\end{equation}
The iteration schemes \tagref{2.18} and \tagref{2.19} yield for $m\ge0$
\begin{equation}
Q_m^L(q) = \sum_{j\ge1} (-1)^{j-1} q^{\frac{j(3j-1)}{2} + mj}
\qbin{2L-m-\au(n,j)}{L+\lfloor\tfrac{n+1}{n}j \rfloor},
\tag{2.20}
\end{equation}
and
\begin{equation}
Q_m^L(q) = \sum_{j\ge1} (-1)^{j-1} q^{\frac{j(3j-1)}{2} + mj}
\qbin{2L+1-m+\au(n,j)}{L+1-m+\lfloor-\tfrac{n+1}{n}j \rfloor},
\tag{2.21}
\end{equation}
respectively.
If we employ \tagref{2.3}, \tagref{2.20} with $m=0$
and \tagref{2.21} with $m=1$ we obtain
\tagref{1.30} with $\sigma=0$. On the other hand, \tagref{2.3}, \tagref{2.20} 
with $m=1$ and $L\to L+1$ together with
\tagref{2.21} with $m=0$ give       
\tagref{1.30} with $\sigma=1$.
To prove \tagref{1.30} with $\sigma=-1$ we note that \tagref{2.7} and 
\tagref{2.3} with $m=-1$ can be combined to give
\begin{equation}
\delta_{m,-1} +
Q_m^L(q) + q^{m+1} Q_{m+3}^{L+1}(q)
= 
q^{m+1} \qbin{2L-m}{L+1},\quad m\ge -1.
\tag{2.22}
\end{equation}
Therefore \tagref{2.20} can be slightly generalized as
\begin{equation}
\delta_{m,-1} + Q_m^L(q) = \sum_{j\ge1} (-1)^{j-1} q^{\frac{j(3j-1)}{2} + mj}
\qbin{2L-m-\au(n,j)}{L+\lfloor\tfrac{n+1}{n}j \rfloor},\quad m\ge -1.
\tag{2.23}
\end{equation}

Next, equation \tagref{2.3} with $m=-1$ and $L\to L-1$ becomes
\begin{equation}
1+ Q_{-1}^{L-1}(q) + Q_{2}^{L}(q) = \qbin{2L-1}{L-1}.
\tag{2.24}
\end{equation}
The last equation together with \tagref{2.21} with $m=2$
and \tagref{2.23} with $m=-1$ and $L\to L-1$ gives
\tagref{1.30} with $\sigma=-1$, as desired.
%%% END INSERT II%%%%%%%%%%%%%%%%%%%%%%%%%%%%%%%%%%%%%%%%%%%%%%%%%%%%%%%%%%%

We now move on to generalize \tagref{1.11}. To this end we define
$P_m^L(q)$ as
\begin{equation}
P_m^L(q) = \sum_{n\ge1} p_m^L(n) q^n,
\tag{2.25}
\end{equation}
where $p_m^L(n)$ is the number of partitions of $n$ with
largest part $\le L$ and rank $m$. Obviously,
\begin{equation}
P_m^L(q) = Q_m^L(q) - Q_{m+1}^L(q).
\tag{2.26}
\end{equation}
So using \tagref{2.10}, \tagref{2.11} and \tagref{2.17}
we obtain
\begin{equation}
P_m^L(q) = 
\sum_{j\ge1} (-1)^{j-1} q^{\frac{j(3j-1)}{2} + mj}
\qbin{2L-m}{L+\lfloor\tfrac{3j}{2} \rfloor}
-\sum_{j\ge1} (-1)^{j-1} q^{\frac{j(3j+1)}{2} + mj}
\qbin{2L-m}{L-\lfloor-\tfrac{3j}{2} \rfloor},
\tag{2.27}
\end{equation}
provided $m\ge0$.
Using the obvious conjugation symmetry
\begin{equation}
P_{-|m|}^L(q) = P_{|m|}^{L+|m|}(q)
\tag{2.28}
\end{equation}
it is straightforward to extend \tagref{2.27} to negative $m$.
This way we obtain the following polynomial analogue of
\tagref{1.11}
\begin{align}
P_m^L(q) &= 
\sum_{j\ge1} (-1)^{j-1} q^{\frac{j(3j-1)}{2} + mj}
\qbin{2L-m}{L+\sign(m)\lfloor\tfrac{3j}{2} \rfloor}\tag{2.29}\\
&\quad -\sum_{j\ge1} (-1)^{j-1} q^{\frac{j(3j+1)}{2} + mj}
\qbin{2L-m}{L-\sign(m)\lfloor-\tfrac{3j}{2} \rfloor},
\nonumber  
\end{align}    
where
$$
\sign(m) = 
\begin{cases}
  1, &\mbox{if $m\ge0$}, \\
  -1, &\mbox{otherwise.}       
\end{cases}
$$
%%%%%END SECTION 2

\setcounter{equation}{0}
\section{Partitions with prescribed cranks}
\label{sec:3}
%%%%%%SECTION 3

Dyson \cite{D1} conjectured that the generating function for the
crank should have a form similar to \tagref{1.11}, and it does
as can be seen from the following formula
\begin{equation}
\Chat_k(q) = 
\frac{1}{(q)_\infty} 
\sum_{j\ge1} (-1)^{j-1} q^{T_{j-1} + j|k|}(1-q^j)
+ q (\delta_{k,0} - \delta_{k,1}),
\tag{3.1}
\end{equation}
where
\begin{equation}
T_j = \frac{j(j+1)}{2},
\tag{3.2}
\end{equation}
and
\begin{equation}
\Chat_k(q) = \sum_{n\ge0} \chat_k(n) q^n,
\tag{3.3}
\end{equation}
with $\chat_k(n)$ denoting the number of partitions of $n$ with
crank $k$. 
In \tagref{3.3} we adopt the convention that $\chat_k(0)=\delta_{k,0}$.
Formula \tagref{3.1} is a consequence of
Theorem (7.19) in \cite{Ga1} and Theorem 1 in \cite{AG}.

To explain \tagref{3.1} in a combinatorial fashion Dyson \cite{D4}
introduced the concept of the rank-set $R(\pi)$ of
a partition
$\pi = p_1 + p_2 + p_3 + \cdots$ with
parts $p_1\ge p_2\ge p_3\ge\cdots$.
$R(\pi)$ is defined as
\begin{equation}
R(\pi) = [j - p_{j+1}, j=0,1,2,\dots].
\tag{3.4}
\end{equation}

To prove \tagref{3.1} Dyson first established that
\begin{equation}
C_k(q) = G_k(q) + q \delta_{k,0},
\tag{3.5}
\end{equation}
where
\begin{align}
C_k(q) = \sum_{n\ge0} c_k(n) q^n, \tag{3.6}\\
G_k(q) = \sum_{n\ge0} g_k(n) q^n, 
\tag{3.7}
\end{align}
with $c_k(n)$, $g_k(n)$ denoting the number of partitions of $n$
with crank $\le k$ and $k$ in the rank-set of these partitions,
respectively. 
In \tagref{3.6}--\tagref{3.7} we use the convention that
$c_k(0)=g_k(0)=1$ if $k\ge0$ and $0$, otherwise.
He then showed that
\begin{equation}
G_{-k}(q) + G_{k-1}(q) = \frac{1}{(q)_\infty},
\tag{3.8}
\end{equation}
and
\begin{equation}
G_{k}(q) + q^{k+1} G_{k+1}(q) = \frac{1}{(q)_\infty},\quad k\ge-1.
\tag{3.9}
\end{equation}

Iteration of \tagref{3.9} yields
\begin{equation}
G_k(q) = \frac{1}{(q)_\infty}
\sum_{j\ge0} (-1)^j q^{T_j + kj}, \quad k\ge-1.
\tag{3.10}
\end{equation}
Now \tagref{3.10}, \tagref{3.5} and the obvious relation
\begin{equation}
\Chat_k(q) = C_k(q) - C_{k-1}(q)
\tag{3.11}
\end{equation}
together imply that
\begin{equation}
\Chat_k(q) = \frac{1}{(q)_\infty}
\sum_{j\ge0} (-1)^{j-1} \left(q^{T_{j-1} + kj} - q^{T_{j} + kj}\right)
+ q(\delta_{k,0} - \delta_{k,1}),\quad k\ge0,
\tag{3.12}
\end{equation}
which is \tagref{3.1} with $k\ge0$. To extend \tagref{3.12} to negative $k$,
we observe that \tagref{3.8} implies that
\begin{equation}
\Ghat_{-k}(q) = \Ghat_k(q),
\tag{3.13}
\end{equation}
where
\begin{equation}
\Ghat_k(q) = G_k(q) - G_{k-1}(q).
\tag{3.14}
\end{equation}
From \tagref{3.5} we deduce that
\begin{equation}
\Chat_k(q) = \Ghat_k(q) + q(\delta_{k,0} - \delta_{k,1}).
\tag{3.15}
\end{equation}
%%If we now replace $k$ by $-k$ in \tagref{3.15} with $k>0$ and use
%%\tagref{3.13} we obtain
%%\begin{equation}
%%\Chat_{-k}(q) = \Ghat_{-k}(q) = \Ghat_k(q),\quad k>0.
%%\tag{3.16}
%%\end{equation}
%%The last equation together with \tagref{3.10}, \tagref{3.14} gives
%%\tagref{3.1} for $k < 0$.
If we now replace $k$ by $-k$ in \tagref{3.15} with $k\ge0$ and use
\tagref{3.13} we obtain
\begin{equation}
\Chat_{-k}(q) = \Ghat_{-k}(q) +q \delta_{k,0}= \Ghat_k(q)+q\delta_{k,0},
\quad k\ge0.
\tag{3.16}
\end{equation}
This equation together with \tagref{3.10}, \tagref{3.14} gives
\tagref{3.1} for $k < 0$.
In addition, using \tagref{3.15} we see that \tagref{3.16} implies
that
\begin{equation}
\Chat_{-k}(q) = \Chat_k(q) + q \delta_{k,1}, \quad k\ge0.
\tag{3.17}
\end{equation}
In the appendix we give a direct combinatorial proof of \tagref{3.17}
without using \tagref{3.15}.

To deal with \tagref{3.8}, \tagref{3.9} Dyson introduced a simple graphical
tool to determine whether or not $k\in R(\pi)$. To explain it we follow
Dyson \cite{D4} and define the boundary of the Ferrers graph of $\pi$
as the infinite zig-zag line consisting of vertical and
horizontal segments each of unit length (see Fig.\ 3).

\begin{figure}[h!]
\begin{center}
\setlength{\unitlength}{0.00083333in}
\begingroup\makeatletter\ifx\SetFigFont\undefined%
\gdef\SetFigFont#1#2#3#4#5{%
  \reset@font\fontsize{#1}{#2pt}%
  \fontfamily{#3}\fontseries{#4}\fontshape{#5}%
  \selectfont}%
\fi\endgroup%
{\renewcommand{\dashlinestretch}{30}
\begin{picture}(3242,2088)(0,-10)
\path(150,1833)(450,1833)(450,1533)
	(150,1533)(150,1833)
\path(450,1833)(750,1833)(750,1533)
	(450,1533)(450,1833)
\path(750,1833)(1050,1833)(1050,1533)
	(750,1533)(750,1833)
\path(1050,1833)(1350,1833)(1350,1533)
	(1050,1533)(1050,1833)
\path(1350,1833)(1650,1833)(1650,1533)
	(1350,1533)(1350,1833)
\path(150,1533)(450,1533)(450,1233)
	(150,1233)(150,1533)
\path(450,1533)(750,1533)(750,1233)
	(450,1233)(450,1533)
\path(150,1233)(450,1233)(450,933)
	(150,933)(150,1233)
\thicklines
\path(2700,1833)(1650,1833)(1650,1533)
	(750,1533)(750,1233)(450,1233)
	(450,933)(150,933)(150,333)
\drawline(2775,1833)(2800,1833)
\drawline(2850,1833)(2875,1833)
\drawline(2925,1833)(2950,1833)
\drawline(3000,1833)(3025,1833)
\drawline(150,250)(150,225)
\drawline(150,175)(150,150)
\drawline(150,100)(150,75)
\drawline(150,25)(150,0)
\put(0,108){$y$}
\put(3150,1983){$x$}
\end{picture}
}
\end{center}
\mycap{3}{\enskip Graph of $\pi=5+2+1$, the boundary $B(\pi)$ is indicated
by the thick line.}
\end{figure}

Next, we draw two $45^{\mbox{o}}$ lines, namely
$$
y=k+x, \qquad y = 1 + k + x,
$$
as shown in Fig.\ 4 and Fig.\ 5.

\begin{figure}[h!]
\begin{center}
\setlength{\unitlength}{0.00083333in}
\begingroup\makeatletter\ifx\SetFigFont\undefined%
\gdef\SetFigFont#1#2#3#4#5{%
  \reset@font\fontsize{#1}{#2pt}%
  \fontfamily{#3}\fontseries{#4}\fontshape{#5}%
  \selectfont}%
\fi\endgroup%
{\renewcommand{\dashlinestretch}{30}
\begin{picture}(2417,2020)(0,-10)
\thicklines
\drawline(225,275)(225,300)
\drawline(225,200)(225,225)
\drawline(225,125)(225,150)
\drawline(225,50)(225,75)
\drawline(2100,1765)(2125,1765)
\drawline(2175,1765)(2200,1765)
\drawline(2250,1765)(2275,1765)
\drawline(2325,1765)(2350,1765)
\put(-350,865){$1+k$}
\put(0,1165){$k$}
\put(225,865){\blacken\ellipse{50}{50}}
\put(225,865){\ellipse{50}{50}}
\put(225,1165){\blacken\ellipse{50}{50}}
\put(225,1165){\ellipse{50}{50}}
\thinlines
\path(225,1765)(525,1765)(525,1465)
	(225,1465)(225,1765)
\path(525,1765)(825,1765)(825,1465)
	(525,1465)(525,1765)
\path(825,1765)(1125,1765)(1125,1465)
	(825,1465)(825,1765)
\path(225,1465)(525,1465)(525,1165)
	(225,1165)(225,1465)
\path(525,1465)(825,1465)(825,1165)
	(525,1165)(525,1465)
\path(225,1165)(525,1165)(525,865)
	(225,865)(225,1165)
\path(225,865)(525,865)(525,565)
	(225,565)(225,865)
\thicklines
\path(2025,1765)(1125,1765)(1125,1465)
	(825,1465)(825,1165)(525,1165)
	(525,565)(225,565)(225,340)
\path(225,865)(1050,40)
\path(225,1165)(1350,40)
\put(0,40){$y$}
\put(2325,1915){$x$}
\end{picture}
}
\end{center}
\mycap{4}{\enskip Graph of $\pi=3+2+1+1$ with $k=2\in R(\pi)$.}
\end{figure}

\begin{figure}[h!]
\begin{center}
\setlength{\unitlength}{0.00083333in}
\begingroup\makeatletter\ifx\SetFigFont\undefined%
\gdef\SetFigFont#1#2#3#4#5{%
  \reset@font\fontsize{#1}{#2pt}%
  \fontfamily{#3}\fontseries{#4}\fontshape{#5}%
  \selectfont}%
\fi\endgroup%
{\renewcommand{\dashlinestretch}{30}
\begin{picture}(2417,2020)(0,-10)
\thicklines
\drawline(225,265)(225,290)
\drawline(225,190)(225,215)
\drawline(225,115)(225,140)
\drawline(225,40)(225,65)
\drawline(2100,1765)(2125,1765)
\drawline(2175,1765)(2200,1765)
\drawline(2250,1765)(2275,1765)
\drawline(2325,1765)(2350,1765)
\put(225,1465){\blacken\ellipse{50}{50}}
\put(225,1465){\ellipse{50}{50}}
\put(225,1165){\blacken\ellipse{50}{50}}
\put(225,1165){\ellipse{50}{50}}
\put(-350,1165){$1+k$}
\put(0,1465){$k$}
\thinlines
\path(225,1765)(525,1765)(525,1465)
	(225,1465)(225,1765)
\path(525,1765)(825,1765)(825,1465)
	(525,1465)(525,1765)
\path(225,1465)(525,1465)(525,1165)
	(225,1165)(225,1465)
\thicklines
\path(2025,1765)(825,1765)(825,1465)
	(525,1465)(525,1165)(225,1165)(225,340)
\path(225,1465)(1050,640)
\path(225,1165)(975,415)
\put(0,40){$y$}
\put(2325,1915){$x$}
\end{picture}
}
\end{center}
\mycap{5}{\enskip Graph of $\pi=2+1$ with $k=1\not\in R(\pi)$.}
\end{figure}

Let $BS_k(\pi)$ denote the segment of $B(\pi)$ lying in the strip
$$
k+x \le y \le k+1+x
$$
determined by these two lines. 
Now if $BS_k(\pi)$ is vertical, then $k\in R(\pi)$, otherwise $k\not\in R(\pi)$.
Using this criterion it is easy to verify that $\nu(\pi)\ne1+k$,
whenever $k\in R(\pi)$.

We are now ready to prove \tagref{3.8}.
First, it is obvious that any given partition $\pi$ counted by 
$\frac{1}{(q)_\infty}$ has either $-k\in R(\pi)$ or $-k\not\in R(\pi)$. 
In the first case, $\pi$ is counted by $G_{-k}(q)$
in \tagref{3.8}.
In the second case, $BS_{-k}(\pi)$ is a horizontal segment, and so if we
conjugate $\pi$ to get $\pi^{*}$, then it is clear that
$BS_{k-1}(\pi^{*})$ is vertical. Therefore, $k-1\in R(\pi^{*})$
and consequently $\pi^{*}$ is counted by $G_{k-1}(q)$ in \tagref{3.8}.

To prove \tagref{3.9}, we remove the row containing the segment $BS_k(\pi)$
from some given partition $\pi$ counted by $G_k(q)$. Next,
we insert a vertical column of height $j+k$ to the right of the
rectangle $[j,j+k]$, where $j$ is the length of the row removed.
This procedure is illustrated in Fig.\ 6.

\begin{figure}[h!]
\begin{center}
\setlength{\unitlength}{0.00071667in}
\begingroup\makeatletter\ifx\SetFigFont\undefined%
\gdef\SetFigFont#1#2#3#4#5{%
  \reset@font\fontsize{#1}{#2pt}%
  \fontfamily{#3}\fontseries{#4}\fontshape{#5}%
  \selectfont}%
\fi\endgroup%
{\renewcommand{\dashlinestretch}{30}
\begin{picture}(6912,3714)(0,-10)
\texture{8101010 10000000 444444 44000000 11101 11000000 444444 44000000 
	101010 10000000 444444 44000000 10101 1000000 444444 44000000 
	101010 10000000 444444 44000000 11101 11000000 444444 44000000 
	101010 10000000 444444 44000000 10101 1000000 444444 44000000 }
\shade\path(675,1812)(1875,1812)(1875,1512)
	(675,1512)(675,1812)
\path(675,1812)(1875,1812)(1875,1512)
	(675,1512)(675,1812)
\shade\path(5400,3687)(5700,3687)(5700,1887)
	(5400,1887)(5400,3687)
\path(5400,3687)(5700,3687)(5700,1887)
	(5400,1887)(5400,3687)
\put(6150,3162){\ellipse{336}{336}}
\put(6100,3087){\makebox(0,0)[lb]{\smash{$2$}}}
\put(675,3012){\ellipse{50}{50}}
\put(675,3012){\blacken\ellipse{50}{50}}
%%%
%%%
\put(675,2712){\ellipse{50}{50}}
\put(675,2712){\blacken\ellipse{50}{50}}
\put(1050,1137){\ellipse{336}{336}}
\put(2325,3087){\ellipse{336}{336}}
\put(4200,3087){\ellipse{50}{50}}
\put(4200,3087){\blacken\ellipse{50}{50}}
%%\put(4200,2787){\ellipse{50}{50}}
%%\put(4200,2787){\blacken\ellipse{50}{50}}
\put(4200,3387){\ellipse{50}{50}}
\put(4200,3387){\blacken\ellipse{50}{50}}
\put(4575,1512){\ellipse{336}{336}}
\path(675,3612)(1875,3612)(1875,1812)
	(675,1812)(675,3612)
\path(2325,1362)(675,3012)
\path(2325,1062)(675,2712)
\path(1350,1662)(1875,1662)
\path(1755.000,1632.000)(1875.000,1662.000)(1755.000,1692.000)
\path(1200,1662)(675,1662)
\path(795.000,1692.000)(675.000,1662.000)(795.000,1632.000)
\path(1905.000,3492.000)(1875.000,3612.000)(1845.000,3492.000)
\path(1875,3612)(1875,1812)
\path(1845.000,1932.000)(1875.000,1812.000)(1905.000,1932.000)
\path(1875,1512)(1875,1212)(1575,1212)
	(1575,912)(1275,912)(1275,312)
	(975,312)(975,12)(675,12)(675,1512)
\path(1875,1812)(2175,1812)(2175,2112)
	(2475,2112)(2475,2712)(2775,2712)
	(2775,3012)(3075,3012)(3075,3612)(1875,3612)
\path(4200,3687)(5400,3687)(5400,1887)
	(4200,1887)(4200,3687)
\path(5850,1437)(4200,3087)
%%\path(5850,1137)(4200,2787)
\path(5850,1737)(4200,3387)
\path(5700,1887)(6000,1887)(6000,2187)
	(6300,2187)(6300,2787)(6600,2787)
	(6600,3087)(6900,3087)(6900,3687)(5700,3687)
\path(5400,1887)(5400,1587)(5100,1587)
	(5100,1287)(4800,1287)(4800,687)
	(4500,687)(4500,387)(4200,387)(4200,1887)
\thicklines
\path(2475,1962)(3375,1962)
\blacken\thinlines
\path(3255.000,1872.000)(3375.000,1962.000)(3255.000,2052.000)(3291.000,1962.000)(3255.000,1872.000)
\blacken\path(4320.000,1917.000)(4200.000,1887.000)(4320.000,1857.000)(4284.000,1887.000)(4320.000,1917.000)
\path(4200,1887)(5400,1887)
\blacken\path(5280.000,1857.000)(5400.000,1887.000)(5280.000,1917.000)(5316.000,1887.000)(5280.000,1857.000)
\blacken\path(5430.000,3567.000)(5400.000,3687.000)(5370.000,3567.000)(5400.000,3603.000)(5430.000,3567.000)
\path(5400,3687)(5400,1887)
\blacken\path(5370.000,2007.000)(5400.000,1887.000)(5430.000,2007.000)(5400.000,1971.000)(5370.000,2007.000)
\put(1225,1587){\makebox(0,0)[lb]{\smash{$j$}}}
\put(1350,2562){\makebox(0,0)[lb]{\smash{$j+k$}}}
\put(2275,3012){\makebox(0,0)[lb]{\smash{$2$}}}
\put(1000,1062){\makebox(0,0)[lb]{\smash{$1$}}}
\put(0,1812){\makebox(0,0)[lb]{\smash{$\pi$:}}}
\put(425,3012){\makebox(0,0)[lb]{\smash{$k$}}}
\put(0,2637){\makebox(0,0)[lb]{\smash{$1+k$}}}
%%\put(4875,2637){\makebox(0,0)[lb]{\smash{$j+k$}}}
\put(4875,2937){\makebox(0,0)[lb]{\smash{$j+k$}}}
\put(3525,1887){\makebox(0,0)[lb]{\smash{$\pi'$:}}}
\put(3950,3087){\makebox(0,0)[lb]{\smash{$k$}}}
%%\put(3600,2712){\makebox(0,0)[lb]{\smash{$1+k$}}}
%%\put(3600,3462){\makebox(0,0)[lb]{\smash{$-1+k$}}} 
%%%\put(3400,3462){\makebox(0,0)[lb]{\smash{$-1+k$}}}
\put(3400,3387){\makebox(0,0)[lb]{\smash{$-1+k$}}}
\put(4525,1437){\makebox(0,0)[lb]{\smash{$1$}}}
\put(4725,1962){\makebox(0,0)[lb]{\smash{$j$}}}
\end{picture}
}
\end{center}
\mycap{6}{\enskip The transformation $\pi\to\pi'$ used in the proof of
\tagref{3.9} ($k\ge0$).}
\end{figure}

Let us call the resulting partition $\pi'$. It is easy to see that
$$
|\pi'| = |\pi| + k,
$$
and, because $BS_{-1+k}(\pi')$ is a horizontal segment,
$$
k-1 \not\in R(\pi').
$$
Since the map $\pi\to\pi'$ is reversible, we immediately infer that
\begin{equation}
g_k(n) = p(n+k) - g_{k-1}(n + k), \qquad k\ge0,
\tag{3.18}
\end{equation}
where $n=|\pi|$. The last equation can be easily transformed in \tagref{3.9}.

In \cite{D4}, Dyson proves \tagref{3.5} first by mapping partitions 
$\pi$ with $k\in R(\pi)$ onto
certain vector partitions introduced in \cite{Ga1}, and then 
mapping these vector partitions onto ordinary partitions with crank $\le k$.
This approach involved ten separate cases. Here, we choose to
prove \tagref{3.5} directly, without any reference to vector partitions.
Our analysis requires consideration of only three separate cases, as we
now explain.

\underbar{Case 1.} Here we consider partitions $\pi$ with $k\in R(\pi)$
and $\nu(\pi)\ge k+2$. This case is illustrated in Fig.\ 7.

\begin{figure}[h!]
\begin{center}
\setlength{\unitlength}{0.00075000in}
\begingroup\makeatletter\ifx\SetFigFont\undefined%
\gdef\SetFigFont#1#2#3#4#5{%
  \reset@font\fontsize{#1}{#2pt}%
  \fontfamily{#3}\fontseries{#4}\fontshape{#5}%
  \selectfont}%
\fi\endgroup%
{\renewcommand{\dashlinestretch}{30}
\begin{picture}(6612,6402)(0,-10)
\texture{8101010 10000000 444444 44000000 11101 11000000 444444 44000000 
	101010 10000000 444444 44000000 10101 1000000 444444 44000000 
	101010 10000000 444444 44000000 11101 11000000 444444 44000000 
	101010 10000000 444444 44000000 10101 1000000 444444 44000000 }
\shade\path(675,4293)(1875,4293)(1875,3993)
	(675,3993)(675,4293)
\path(675,4293)(1875,4293)(1875,3993)
	(675,3993)(675,4293)
\shade\path(4200,1668)(4500,1668)(4500,2868)
	(4200,2868)(4200,1668)
\path(4200,1668)(4500,1668)(4500,2868)
	(4200,2868)(4200,1668)
\put(712.500,5080.500){\arc{237.171}{0.3218}{1.8925}}
\put(4237.500,5455.500){\arc{237.171}{0.3218}{1.8925}}
\put(675,5493){\blacken\ellipse{50}{50}}
\put(675,5493){\ellipse{50}{50}}
\put(675,5193){\blacken\ellipse{50}{50}}
\put(675,5193){\ellipse{50}{50}}
\put(4200,5568){\blacken\ellipse{50}{50}}
\put(4200,5568){\ellipse{50}{50}}
\put(1050,3618){\ellipse{336}{336}}
\put(2325,5568){\ellipse{336}{336}}
\put(4575,3993){\ellipse{336}{336}}
\put(5925,5643){\ellipse{336}{336}}
\path(675,6093)(1875,6093)(1875,4293)
	(675,4293)(675,6093)
\path(2325,3843)(675,5493)
\path(2325,3543)(675,5193)
\path(1875,4293)(2175,4293)(2175,4593)
	(2475,4593)(2475,5193)(2775,5193)
	(2775,5493)(3075,5493)(3075,6093)(1875,6093)
\path(4200,6168)(5400,6168)(5400,4368)
	(4200,4368)(4200,6168)
\path(5850,3918)(4200,5568)
\thicklines
\path(2475,4443)(3375,4443)
\blacken\thinlines
\path(3255.000,4353.000)(3375.000,4443.000)(3255.000,4533.000)(3291.000,4443.000)(3255.000,4353.000)
\blacken\path(5430.000,6048.000)(5400.000,6168.000)(5370.000,6048.000)(5400.000,6084.000)(5430.000,6048.000)
\path(5400,6168)(5400,4368)
\blacken\path(5370.000,4488.000)(5400.000,4368.000)(5430.000,4488.000)(5400.000,4452.000)(5370.000,4488.000)
\path(5400,4368)(5700,4368)(5700,4668)
	(6000,4668)(6000,5268)(6300,5268)
	(6300,5568)(6600,5568)(6600,6168)(5400,6168)
\path(4200,4368)(5400,4368)
\path(4200,4368)(4200,2868)(4500,2868)
	(4500,3168)(4800,3168)(4800,3768)
	(5100,3768)(5100,4068)(5400,4068)(5400,4368)
\path(675,3993)(675,2493)(975,2493)
	(975,2793)(1275,2793)(1275,3393)
	(1575,3393)(1575,3693)(1875,3693)(1875,3993)
\blacken\path(4680.000,2748.000)(4650.000,2868.000)(4620.000,2748.000)(4650.000,2784.000)(4680.000,2748.000)
\path(4650,2868)(4650,2418)
\blacken\path(4020.000,1788.000)(4050.000,1668.000)(4080.000,1788.000)(4050.000,1752.000)(4020.000,1788.000)
\path(4050,1668)(4050,2268)
\dashline{60.000}(4500,3168)(3975,3168)
\blacken\path(795.000,6123.000)(675.000,6093.000)(795.000,6063.000)(759.000,6093.000)(795.000,6123.000)
\path(675,6093)(1875,6093)
\blacken\path(1755.000,6063.000)(1875.000,6093.000)(1755.000,6123.000)(1791.000,6093.000)(1755.000,6063.000)
\blacken\path(1905.000,5973.000)(1875.000,6093.000)(1845.000,5973.000)(1875.000,6009.000)(1905.000,5973.000)
\path(1875,6093)(1875,4293)
\blacken\path(1845.000,4413.000)(1875.000,4293.000)(1905.000,4413.000)(1875.000,4377.000)(1845.000,4413.000)
\blacken\path(4080.000,3048.000)(4050.000,3168.000)(4020.000,3048.000)(4050.000,3084.000)(4080.000,3048.000)
\path(4050,3168)(4050,2493)
\blacken\path(4620.000,1788.000)(4650.000,1668.000)(4680.000,1788.000)(4650.000,1752.000)(4620.000,1788.000)
\path(4650,1668)(4650,2193)
\blacken\path(4320.000,6198.000)(4200.000,6168.000)(4320.000,6138.000)(4284.000,6168.000)(4320.000,6198.000)
\path(4200,6168)(5400,6168)
\blacken\path(5280.000,6138.000)(5400.000,6168.000)(5280.000,6198.000)(5316.000,6168.000)(5280.000,6138.000)
\put(1350,5043){\makebox(0,0)[lb]{\smash{$j+k$}}}
\put(2275,5493){\makebox(0,0)[lb]{\smash{$2$}}}
\put(1000,3543){\makebox(0,0)[lb]{\smash{$1$}}}
\put(0,4293){\makebox(0,0)[lb]{\smash{$\pi$:}}}
\put(475,5493){\makebox(0,0)[lb]{\smash{$k$}}}
\put(75,5118){\makebox(0,0)[lb]{\smash{$1+k$}}}
\put(4875,5118){\makebox(0,0)[lb]{\smash{$j+k$}}}
\put(3525,4368){\makebox(0,0)[lb]{\smash{$\pi'$:}}}
\put(4000,5568){\makebox(0,0)[lb]{\smash{$k$}}}
\put(4525,3918){\makebox(0,0)[lb]{\smash{$1$}}}
\put(1200,6168){\makebox(0,0)[lb]{\smash{$j$}}}
\put(750,4743){\makebox(0,0)[lb]{\smash{$45^{\mbox{o}}$}}}
\put(5875,5568){\makebox(0,0)[lb]{\smash{$2$}}}
\put(4800,6243){\makebox(0,0)[lb]{\smash{$j$}}}
\put(4650,2268){\makebox(0,0)[lb]{\smash{$j$}}}
\put(3975,2343){\makebox(0,0)[lb]{\smash{$\mu$}}}
\put(4275,5118){\makebox(0,0)[lb]{\smash{$45^{\mbox{o}}$}}}
\end{picture}
}
\end{center}
\vskip -60pt
\textsc{Fig 7}.\ \parbox[t]{5in}{Map $\pi\to\pi'$ from partitions $\pi$ with
$k\in R(\pi)$, $\nu(\pi)\ge k+2$ to partitions $\pi'$ with $\crank(\pi')\le k$,
$\mu(\pi')>0$, $\nu(\pi')\ge k+2$.}
\end{figure}

We now remove the row bounded by the vertical segment $BS_k(\pi)$
and then add a vertical column representing $j$ ones to the resulting
graph, where $j>0$ is the length of the row removed. We call this
last partition $\pi'$. It is easy to see that
$$
\nu(\pi')\ge k+2,\quad
\mu(\pi')\ge j >0,\quad\mbox{and}\quad \nutwid(\pi') \le j+k,
$$
where $\mu$ and $\nutwid$ were defined in \tagref{1.6}.
Clearly, $\crank(\pi')=\nutwid(\pi') - \mu(\pi')\le k$.
Perhaps, it is not immediately obvious that the map $\pi\to\pi'$
is reversible. To see that it is, we consider partitions $\pi'$
with $\crank(\pi')\le k$, $\mu(\pi')>0$ and $\nu(\pi')\ge k+2$.
Next we define $j$ to be the $x$-coordinate of the intersection
point of the line $y=x+k$ and the boundary $B(\pi')$.
Since $\nu(\pi')\ge k+2$, $j$ is positive. Moreover, $j\le\mu$
because otherwise $\crank(\pi')$ would be $>k$. So we can remove from
$\pi'$ a vertical column of length $j$ representing ones and place it as
a row of length $j$ right underneath the $[j,j+k]$ rectangle.
This way we obtain $\pi$ with $k\in R(\pi)$, $\nu(\pi)\ge k+2$.

\underbar{Case 2.} Here we consider partitions $\pi$ with $\nu(\pi)\le k$
and unique largest part $\lambda(\pi)$. In this case the
segment $BS_k(\pi)$ is necessarily vertical, implying that
$k\in R(\pi)$. We now transform $\pi$ into $\pi'$ as follows.
If $|\pi|>1$, then we add a part of size $1$ to $\pi$
and subtract $1$ from $\lambda(\pi)$, giving
$\lambda(\pi')=\lambda(\pi)-1$, $\nu(\pi')=\nu(\pi)+1$, $\mu(\pi')>0$.
If $|\pi|=1$, then we define $\pi'=\pi$. It is obvious that the
map $\pi\to\pi'$ is reversible and that $\crank(\pi')\le k-1$, $\mu(\pi')>0$,
$\nu(\pi')\le k+1$.

\underbar{Case 3.} Here we consider partitions $\pi$ with $k\ge2$, 
$\nu(\pi)\le k$ and the largest part $\lambda(\pi)$ is repeated.
Once again, it is clear that $k\in R(\pi)$. We now conjugate $\pi$
to get $\pi'=\pi^{*}$. Since the smallest part of $\pi'$ is
at least $2$, we have $\mu(\pi')=0$, and 
$\crank(\pi')=\lambda(\pi') =\nu(\pi)\le k$.

We now recall that $\nu(\pi)\ne k+1$ whenever $k\in R(\pi)$.
Thus the three cases above are exhaustive. Hence, \tagref{3.5} holds
for $k>0$.

If $k<0$ there is no need to consider cases 2 and 3, because there
are no partitions with a negative number of parts. In addition,
case 1 requires no modification. Hence, \tagref{3.5}
is valid in this case $k<0$, as well.

If $k=0$, then there is no need to consider case 3. Once again, case 1
requires no modification. However, in case 2 the map 
$\pi\to\pi'$ is not bijective. To see this, we note that the set
of partitions $\pi$ with $\nu(\pi)\le0$ is empty, but the set of partitions
$\pi'$ with $\crank(\pi')\le-1$, $\mu(\pi')>0$,
$\nu(\pi')\le 1$ consists of the single partition $\pi'$ with $|\pi'|=1$,
$\nu(\pi')=1$ and $\crank(\pi')=-1$.
Thus,
\begin{equation}
c_0(n) = g_0(n) + \delta_{n,1}, \qquad n\ge0.
\tag{3.19}
\end{equation}
The last equation can be easily transformed into \tagref{3.5} with $k=0$.

%%%%%END SECTION 3

\setcounter{equation}{0}
\section{Partitions with bounds on the largest part and the crank}
\label{sec:4}
%%%%%%SECTION 4

Let $C_k^L(q)$, $\Chat_k^L(q)$, $G_k^L(q)$ denote the generating functions
for partitions with crank $\le k$ and  largest part $\le L$, with crank $k$
and largest part $\le L$, with $k$ in the rank-set and  largest part $\le L$,
respectively.  In this section we will establish the following bounded
analogues of \tagref{3.5} and \tagref{3.9}:
\begin{align}    
&C_k^L(q) = G_k^L(q) + \frac{1-q}{(q)_k} + (q-1)\qbin{L+k}{k} +q \delta_{k,0},
\tag{4.1}\\
&G_{k}^L(q) + q^{k+1} G_{k+1}^{L-1}(q) = \frac{1}{(q)_L},
\tag{4.2}
\end{align}    
where for the sake of simplicity here (and throughout this section) we assume 
that $0\le k \le L$,
$L\ne0$, unless otherwise stated.

The proof of \tagref{4.2} is essentially the same as that of \tagref{3.9}.
Iterating \tagref{4.2} we derive
\begin{equation}
G_k^L(q) = \sum_{j=0}^L (-1)^j \frac{q^{T_j + kj}}{(q)_{L-j}}.
\tag{4.3}
\end{equation}
To prove \tagref{4.1} we need to follow the three separate cases of
the map $\pi\to\pi'$ we used to prove \tagref{3.5}.

Case 1 requires no modification. In case 2 the map $\pi\to\pi'$
produces partitions $\pi'$ with $\lambda(\pi')=\lambda(\pi)-1\le L-1$
and, therefore, misses partitions $\pi'$ counted by $C_k^L(q)$ in \tagref{4.1}
such that $\lambda(\pi')=L$, $\mu(\pi')>0$ and $\nu(\pi')\le k+1$,
and when $k=0$ this map also misses 
the partition $\pi'=1$, as discussed earlier. In other words, the
correction term needed in this case is
\begin{equation}
\CT_2 = q^{1+L}\qbin{L+k-1}{k-1} + q \delta_{k,0}.
\tag{4.4}
\end{equation}
In case 3, the map $\pi\to\pi'$ fails to account for partitions $\pi'$
counted by $C_k^L(q)$ such that $\lambda(\pi')\le k$, $\nu(\pi')>L$,
$\mu(\pi')=0$. The correction term needed in this case is
\begin{equation}
\CT_3 = \left\{ \frac{1-q}{(q)_k} - \left( \qbin{L+k}{k} - q\qbin{L-1+k}{k}
\right)\right\} \theta(k>1), 
\tag{4.5}
\end{equation}
where
\begin{equation}
\theta(\mbox{statement}) =
\begin{cases}
  1, &\mbox{if statement is true}, \\
  0, &\mbox{otherwise.}
\end{cases}
\tag{4.6}
\end{equation}
To understand \tagref{4.5} we observe that
$\frac{1-q}{(q)_k}$ is the generating function for partitions without ones
and largest part not exceeding $k$, and
$$
\left( \qbin{L+k}{k} - q\qbin{L-1+k}{k} \right)
$$
is the generating function for partitions $\pitwid$
with 
$\lambda(\pitwid)\le k$, $\nu(\pitwid)\le L$,
$\mu(\pitwid)=0$.
Combining \tagref{4.4} and \tagref{4.5} and using the $q$-binomial recurrence
\tagref{2.6} we get the total correction term
\begin{equation}
T = \CT_2 + \CT_3 = \frac{1-q}{(q)_k} 
+ (q-1)\qbin{L+k}{k} +q \delta_{k,0},
\tag{4.7}
\end{equation}
as desired. Since
\begin{equation}
\Chat_k^L(q) = C_k^L(q) - C_{k-1}^L(q)
\tag{4.8}
\end{equation}
we have for $L\ge k>0$
\begin{align}
\Chat_k^L(q) &= 
 \sum_{j=1}^L (-1)^{j-1} q^{T_{j-1} + kj}\frac{(1-q^j)}{(q)_{L-j}}\tag{4.9}\\
 &\quad
 + \theta(k>1) \frac{1-q}{(q)_k} q^k + (q-1) q^k \qbin{L-1+k}{k},
\nonumber
\end{align}
by using \tagref{4.1} and \tagref{4.3}.

We now derive a very different representation for $\Chat_k^L(q)$ using
\tagref{1.6}. Because the crank is defined in \tagref{1.6}
in a piece-wise fashion we have to treat two separate cases.

\underbar{Case A.} Here we consider partitions $\pio$
with $\crank(\pio)=k>0$, 
$\lambda(\pio)\le L$, and 
$\mu(\pio)>0$. We decompose the graph of some given $\pio$ as shown
in Fig.\ 8 below.

\begin{figure}[h!]
\begin{center}
\setlength{\unitlength}{0.00075000in}
\begingroup\makeatletter\ifx\SetFigFont\undefined%
\gdef\SetFigFont#1#2#3#4#5{%
  \reset@font\fontsize{#1}{#2pt}%
  \fontfamily{#3}\fontseries{#4}\fontshape{#5}%
  \selectfont}%
\fi\endgroup%
{\renewcommand{\dashlinestretch}{30}
\begin{picture}(5298,3339)(0,-10)
\put(612,1662){\ellipse{212}{212}}
\put(3012,2862){\ellipse{212}{212}}
\blacken\path(132.000,3342.000)(12.000,3312.000)(132.000,3282.000)(96.000,3312.000)(132.000,3342.000)
\path(12,3312)(2412,3312)
\blacken\path(2292.000,3282.000)(2412.000,3312.000)(2292.000,3342.000)(2328.000,3312.000)(2292.000,3282.000)
\blacken\path(2442.000,3192.000)(2412.000,3312.000)(2382.000,3192.000)(2412.000,3228.000)(2442.000,3192.000)
\path(2412,3312)(2412,2112)
\blacken\path(2382.000,2232.000)(2412.000,2112.000)(2442.000,2232.000)(2412.000,2196.000)(2382.000,2232.000)
\path(12,3312)(12,2112)(2412,2112)
\path(12,2112)(12,912)(1212,912)
	(1212,1212)(1512,1212)(1512,1512)
	(1812,1512)(1812,1812)(2262,1812)(2262,2112)
\path(12,912)(162,912)(162,762)
	(12,762)(12,912)
\path(12,762)(162,762)(162,612)
	(12,612)(12,762)
\path(12,612)(162,612)(162,462)
	(12,462)(12,612)
\path(12,462)(162,462)(162,312)
	(12,312)(12,462)
\path(12,312)(162,312)(162,162)
	(12,162)(12,312)
\path(12,162)(162,162)(162,12)
	(12,12)(12,162)
\path(2487,2262)(3312,2262)(3312,2412)
	(3462,2412)(3462,2562)(3762,2562)
	(3762,2862)(4212,2862)(4212,3162)
	(4812,3162)(4812,3312)(2412,3312)
\put(837,3087){\makebox(0,0)[lb]{\smash{$\mu + 1$}}}
\put(1862,2637){\makebox(0,0)[lb]{\smash{$\mu + k$}}}
\put(237,387){\makebox(0,0)[lb]{\smash{$\mu > 0$}}}
\put(2412,1887){\makebox(0,0)[lb]{\smash{$\le \mu$}}}
\put(5037,3162){\makebox(0,0)[lb]{\smash{$\le L$}}}
\put(572,1597){\makebox(0,0)[lb]{\smash{{{\SetFigFont{11}{13.2}{\rmdefault}{\mddefault}{\updefault}1}}}}}
\put(2972,2797){\makebox(0,0)[lb]{\smash{{{\SetFigFont{11}{13.2}{\rmdefault}{\mddefault}{\updefault}2}}}}}
\end{picture}
}
\end{center}
%%\vskip -60pt
\begin{center}
\textsc{Fig 8}.\ \parbox[t]{4in}{Decomposition of partition $\pio$
with $\crank(\pio)=k>0$, $L\ge\lambda(\pio)$, $\mu(\pio)>0$.}
\end{center}
\end{figure}

From this decomposition it is clear that the generating function for these
partitions is
\begin{equation}
A(q) = \sum_{\mu=1}^{L-1} q^{\mu} q^{(\mu+1)(\mu+k)}
\frac{1}{(q^2;q)_{\mu-1}} \qbin{L-1+k}{\mu + k}.
\tag{4.10}
\end{equation}

\underbar{Case B.} 
Here we consider partitions $\pio$ without ones with 
$\crank(\pio)=\lambda(\pio)=k$, $2\le k\le L$.
Clearly, the generating function for these partitions is
\begin{equation}
B(q) = \frac{q^k}{(q^2;q)_{k-1}} \theta(k>1).
\tag{4.11}
\end{equation}
Combining \tagref{4.10}, \tagref{4.11} we find that
\begin{align}
\Chat_k^L(q) &= A(q) + B(q) \tag{4.12}\\
&= \frac{q^k(1-q)}{(q)_k} \theta(k>1)
+ \sum_{\mu=1}^{L-1} \frac{q^{(\mu+1)(\mu+k) + \mu}}
{(q^2;q)_{\mu-1}} \qbin{L-1+k}{\mu + k},
\quad 0< k \le L.
\nonumber 
\end{align}
Comparing \tagref{4.9} and \tagref{4.12} we arrive at the following
identity
\begin{equation}
\sum_{j=1}^L (-1)^{j-1} q^{T_{j-1} + kj}\frac{(1-q^j)}{(q)_{L-j}}
= (1-q)
\sum_{\mu=0}^{L-1} \frac{q^{(\mu+1)(\mu+k) + \mu}}
{(q)_{\mu}} \qbin{L-1+k}{\mu + k}.
\tag{4.13}
\end{equation}
Remarkably, this identity is nothing else but a limiting case of
Heine's second transformation of a ${}_2\phi_1$-series \cite{GR1}:
\begin{equation}
\2phi1{a,b}{c}{q,z} = \frac{ \qinf{\frac{c}{b}} \qinf{bz} }
	                   { \qinf{c} \qinf{z} }
			\2phi1{\frac{abz}{c},b}{bz}{q,\frac{c}{b}},
\tag{4.14}
\end{equation}
where
\begin{equation}
\2phi1{a,b}{c}{q,z} = \sum_{n=1}^\infty \frac{ (a)_n (b)_n }
					     { (c)_n (q)_n }
					z^n.
\tag{4.15}
\end{equation}
To see this we rewrite the left side of \tagref{4.13} in
$q$-hypergeometric form as
\begin{equation}
\sum_{j=1}^L (-1)^{j-1} q^{T_{j-1} + kj}\frac{(1-q^j)}{(q)_{L-j}}
= \frac{q^k(1-q)}{(q)_{L-1}}
  \lim_{c\to0} \2phi1{q^2,q^{1-L}}{c}{q,q^{L+k}},
\tag{4.16}
\end{equation}
Here we have used
\begin{equation}
(q)_{L-j} = \frac{ (q)_{L-1}}{ (q^{1-L} )_{j-1}}
		(-1)^{j-1} q^{T_{j-2} - (L-1)(j-1)},
\tag{4.17}
\end{equation}
and
\begin{equation}
\frac{1-q^{1+j}}{1-q} =  \frac{ (q^2)_j }{ (q)_j },
\tag{4.18}
\end{equation}
along with the trivial relation
\begin{equation}
\lim_{c\to0} (c)_n = 1.
\tag{4.19}
\end{equation}
Next we employ \tagref{4.14} with $a=q^2$, $b=q^{1-L}$, $z=q^{L+k}$
together with
\begin{equation}
\frac{ \qinf{q^{1+k}} }{ \qinf{q^{L+k}} } = (q^{1+k})_{L-1},
\tag{4.20}
\end{equation}
and
\begin{equation}
\lim_{\rho\to\infty} (\rho)_i \rho^{-i} = (-1)^i q^{T_{i-1}}
\tag{4.21}
\end{equation}
to derive
\begin{align}    
&\sum_{j=1}^L (-1)^{j-1} q^{T_{j-1} + kj}\frac{(1-q^j)}{(q)_{L-j}}
= \frac{q^k(1-q)}{(q)_{L-1}} (q^{1+k})_{L-1}
  \lim_{c\to0} \2phi1{\frac{q^{3+k}}{c},q^{1-L}}{q^{1+k}}{q,cq^{L-1}}
\tag{4.22}\\
&=
\frac{q^k(1-q)}{(q)_{L-1}} (q^{1+k})_{L-1}
\sum_{i=0}^{L-1} (-1)^i \frac{ (q^{1-L})_i }{ (q)_i (q^{1+k})_i }
			q^{ T_{i-1} + (L+k+2)i }.
\nonumber
\end{align}
Finally, verifying that
\begin{equation}
(-1)^i \frac{ (q^{1+k})_{L-1} }{ (q)_{L-1} } \frac{ (q^{1-L})_i }{ (q^{1+k})_i}
= q^{T_i - Li} \qbin{L-1+k}{i+k}
\tag{4.23}
\end{equation}
we see that
\begin{equation}
\sum_{j=1}^L (-1)^{j-1} q^{T_{j-1} + kj}\frac{(1-q^j)}{(q)_{L-j}}
= (1-q) \sum_{i=0}^{L-1} \frac{ q^{i^2 + (2+k)i + k} }{ (q)_i }
			\qbin{L-1+k}{i+k}.
\tag{4.24}
\end{equation}
This last equation is essentially \tagref{4.13}, as desired.

The $q$-hypergeometric proof of \tagref{4.13} clearly suggests that our
analysis can be extended further to treat partitions $\pi$ with $\crank(\pi)=k$,
$\lambda(\pi)\le L$ and $\nu(\pi) \le M$. However, we will not pursue this
here.

%%%%%END SECTION 4

\setcounter{equation}{0}
\section{A variant of Dyson's transformation \\
and a new proof of Gauss's formula}
\label{sec:5}

%%%%%%SECTION 5

Let $e(n)$ denote the number of partitions of $n$
into distinct odd parts with all other parts being even.
The generating function $E(q)$ for these partitions can be written
in the form of a product as
\begin{equation}
E(q) = \sum_{n\ge0} e(n) q^n = \frac{ (-q;q^2)_\infty }{ (q^2;q^2)_\infty }.
\tag{5.1}
\end{equation}
We will use MacMahon's graphs with modulus $2$ to depict these partitions.
For example, the mod $2$ graph of the partition
$\pi = 7 + 6 + 6 + 5 + 2$ is given in Fig.\ 9.

\begin{figure}[h!]
\begin{center}
\setlength{\unitlength}{0.00083333in}
\begingroup\makeatletter\ifx\SetFigFont\undefined%
\gdef\SetFigFont#1#2#3#4#5{%
  \reset@font\fontsize{#1}{#2pt}%
  \fontfamily{#3}\fontseries{#4}\fontshape{#5}%
  \selectfont}%
\fi\endgroup%
{\renewcommand{\dashlinestretch}{30}
\begin{picture}(4374,1539)(0,-10)
\path(12,312)(312,312)(312,12)
	(12,12)(12,312)
\path(12,612)(312,612)(312,312)
	(12,312)(12,612)
\path(612,612)(912,612)(912,312)
	(612,312)(612,612)
\path(612,912)(912,912)(912,612)
	(612,612)(612,912)
\path(312,612)(612,612)(612,312)
	(312,312)(312,612)
\path(312,912)(612,912)(612,612)
	(312,612)(312,912)
\path(612,1212)(912,1212)(912,912)
	(612,912)(612,1212)
\path(12,912)(312,912)(312,612)
	(12,612)(12,912)
\path(12,1212)(312,1212)(312,912)
	(12,912)(12,1212)
\path(312,1212)(612,1212)(612,912)
	(312,912)(312,1212)
\path(12,1512)(312,1512)(312,1212)
	(12,1212)(12,1512)
\path(312,1512)(612,1512)(612,1212)
	(312,1212)(312,1512)
\path(612,1512)(912,1512)(912,1212)
	(612,1212)(612,1512)
\path(912,1512)(1212,1512)(1212,1212)
	(912,1212)(912,1512)
\path(2262,1512)(2562,1512)(2562,1212)
	(2262,1212)(2262,1512)
\path(2562,1512)(2862,1512)(2862,1212)
	(2562,1212)(2562,1512)
\path(2862,1512)(3162,1512)(3162,1212)
	(2862,1212)(2862,1512)
\path(3162,1512)(3462,1512)(3462,1212)
	(3162,1212)(3162,1512)
\path(2262,1212)(2562,1212)(2562,912)
	(2262,912)(2262,1212)
\path(2562,1212)(2862,1212)(2862,912)
	(2562,912)(2562,1212)
\path(2862,1212)(3162,1212)(3162,912)
	(2862,912)(2862,1212)
\path(3162,1212)(3462,1212)(3462,912)
	(3162,912)(3162,1212)
\path(2262,912)(2562,912)(2562,612)
	(2262,612)(2262,912)
\path(2562,912)(2862,912)(2862,612)
	(2562,612)(2562,912)
\path(2862,912)(3162,912)(3162,612)
	(2862,612)(2862,912)
\path(3162,912)(3462,912)(3462,612)
	(3162,612)(3162,912)
\path(2262,612)(2562,612)(2562,312)
	(2262,312)(2262,612)
\path(2562,612)(2862,612)(2862,312)
	(2562,312)(2562,612)
\path(2862,612)(3162,612)(3162,312)
	(2862,312)(2862,612)
\path(3162,612)(3462,612)(3462,312)
	(3162,312)(3162,612)
\path(3462,1512)(3762,1512)(3762,1212)
	(3462,1212)(3462,1512)
\path(2262,312)(2562,312)(2562,12)
	(2262,12)(2262,312)
\path(2562,312)(2862,312)(2862,12)
	(2562,12)(2562,312)
\path(3462,912)(3762,912)(3762,612)
	(3462,612)(3462,912)
\path(3762,912)(4062,912)(4062,612)
	(3762,612)(3762,912)
\path(3462,1212)(3762,1212)(3762,912)
	(3462,912)(3462,1212)
\path(3762,1212)(4062,1212)(4062,912)
	(3762,912)(3762,1212)
\path(3762,1512)(4062,1512)(4062,1212)
	(3762,1212)(3762,1512)
\path(4062,1512)(4362,1512)(4362,1212)
	(4062,1212)(4062,1512)
\path(3462,612)(3762,612)(3762,312)
	(3462,312)(3462,612)
\put(112,112){\makebox(0,0)[lb]{\smash{{{\SetFigFont{12}{14.4}{\rmdefault}{\mddefault}{\updefault}2}}}}}
\put(112,412){\makebox(0,0)[lb]{\smash{{{\SetFigFont{12}{14.4}{\rmdefault}{\mddefault}{\updefault}2}}}}}
\put(712,412){\makebox(0,0)[lb]{\smash{{{\SetFigFont{12}{14.4}{\rmdefault}{\mddefault}{\updefault}1}}}}}
\put(712,712){\makebox(0,0)[lb]{\smash{{{\SetFigFont{12}{14.4}{\rmdefault}{\mddefault}{\updefault}2}}}}}
\put(412,412){\makebox(0,0)[lb]{\smash{{{\SetFigFont{12}{14.4}{\rmdefault}{\mddefault}{\updefault}2}}}}}
\put(412,712){\makebox(0,0)[lb]{\smash{{{\SetFigFont{12}{14.4}{\rmdefault}{\mddefault}{\updefault}2}}}}}
\put(712,1012){\makebox(0,0)[lb]{\smash{{{\SetFigFont{12}{14.4}{\rmdefault}{\mddefault}{\updefault}2}}}}}
\put(112,712){\makebox(0,0)[lb]{\smash{{{\SetFigFont{12}{14.4}{\rmdefault}{\mddefault}{\updefault}2}}}}}
\put(112,1012){\makebox(0,0)[lb]{\smash{{{\SetFigFont{12}{14.4}{\rmdefault}{\mddefault}{\updefault}2}}}}}
\put(412,1012){\makebox(0,0)[lb]{\smash{{{\SetFigFont{12}{14.4}{\rmdefault}{\mddefault}{\updefault}2}}}}}
\put(112,1312){\makebox(0,0)[lb]{\smash{{{\SetFigFont{12}{14.4}{\rmdefault}{\mddefault}{\updefault}2}}}}}
\put(412,1312){\makebox(0,0)[lb]{\smash{{{\SetFigFont{12}{14.4}{\rmdefault}{\mddefault}{\updefault}2}}}}}
\put(712,1312){\makebox(0,0)[lb]{\smash{{{\SetFigFont{12}{14.4}{\rmdefault}{\mddefault}{\updefault}2}}}}}
\put(1012,1312){\makebox(0,0)[lb]{\smash{{{\SetFigFont{12}{14.4}{\rmdefault}{\mddefault}{\updefault}1}}}}}
\end{picture}
}
\end{center}
%%\vskip -60pt
\textsc{Fig 9}.\ \parbox[t]{5in}{mod $2$ and regular mod $1$ representations
of $\pi = 7 + 6 + 6 + 5 + 2$}
\end{figure}

A nice thing about mod $2$ representations of the partitions counted
by $E(q)$ is that these representations have certain invariance
properties under conjugation. Namely, if we conjugate
the mod $2$ graph of some given partition counted by $E(q)$
we obtain a partition that is also counted by $E(q)$. For instance
if we conjugate the mod $2$ graph of the partition depicted in Fig.\ 9
we get $\pi^{*}=10+8+7+1$ whose mod $2$ graph is given in Fig.\ 10.

\begin{figure}[h!]
\begin{center}
\setlength{\unitlength}{0.00083333in}
\begingroup\makeatletter\ifx\SetFigFont\undefined%
\gdef\SetFigFont#1#2#3#4#5{%
  \reset@font\fontsize{#1}{#2pt}%
  \fontfamily{#3}\fontseries{#4}\fontshape{#5}%
  \selectfont}%
\fi\endgroup%
{\renewcommand{\dashlinestretch}{30}
\begin{picture}(1524,1239)(0,-10)
\path(12,312)(312,312)(312,12)
	(12,12)(12,312)
\path(12,612)(312,612)(312,312)
	(12,312)(12,612)
\path(312,612)(612,612)(612,312)
	(312,312)(312,612)
\path(612,612)(912,612)(912,312)
	(612,312)(612,612)
\path(912,612)(1212,612)(1212,312)
	(912,312)(912,612)
\path(912,912)(1212,912)(1212,612)
	(912,612)(912,912)
\path(612,912)(912,912)(912,612)
	(612,612)(612,912)
\path(312,912)(612,912)(612,612)
	(312,612)(312,912)
\path(12,912)(312,912)(312,612)
	(12,612)(12,912)
\path(12,1212)(312,1212)(312,912)
	(12,912)(12,1212)
\path(312,1212)(612,1212)(612,912)
	(312,912)(312,1212)
\path(612,1212)(912,1212)(912,912)
	(612,912)(612,1212)
\path(912,1212)(1212,1212)(1212,912)
	(912,912)(912,1212)
\path(1212,1212)(1512,1212)(1512,912)
	(1212,912)(1212,1212)
\put(112,112){\makebox(0,0)[lb]{\smash{{{\SetFigFont{12}{14.4}{\rmdefault}{\mddefault}{\updefault}1}}}}}
\put(112,412){\makebox(0,0)[lb]{\smash{{{\SetFigFont{12}{14.4}{\rmdefault}{\mddefault}{\updefault}2}}}}}
\put(412,412){\makebox(0,0)[lb]{\smash{{{\SetFigFont{12}{14.4}{\rmdefault}{\mddefault}{\updefault}2}}}}}
\put(712,412){\makebox(0,0)[lb]{\smash{{{\SetFigFont{12}{14.4}{\rmdefault}{\mddefault}{\updefault}2}}}}}
\put(1012,412){\makebox(0,0)[lb]{\smash{{{\SetFigFont{12}{14.4}{\rmdefault}{\mddefault}{\updefault}1}}}}}
\put(1012,712){\makebox(0,0)[lb]{\smash{{{\SetFigFont{12}{14.4}{\rmdefault}{\mddefault}{\updefault}2}}}}}
\put(712,712){\makebox(0,0)[lb]{\smash{{{\SetFigFont{12}{14.4}{\rmdefault}{\mddefault}{\updefault}2}}}}}
\put(412,712){\makebox(0,0)[lb]{\smash{{{\SetFigFont{12}{14.4}{\rmdefault}{\mddefault}{\updefault}2}}}}}
\put(112,712){\makebox(0,0)[lb]{\smash{{{\SetFigFont{12}{14.4}{\rmdefault}{\mddefault}{\updefault}2}}}}}
\put(112,1012){\makebox(0,0)[lb]{\smash{{{\SetFigFont{12}{14.4}{\rmdefault}{\mddefault}{\updefault}2}}}}}
\put(412,1012){\makebox(0,0)[lb]{\smash{{{\SetFigFont{12}{14.4}{\rmdefault}{\mddefault}{\updefault}2}}}}}
\put(712,1012){\makebox(0,0)[lb]{\smash{{{\SetFigFont{12}{14.4}{\rmdefault}{\mddefault}{\updefault}2}}}}}
\put(1012,1012){\makebox(0,0)[lb]{\smash{{{\SetFigFont{12}{14.4}{\rmdefault}{\mddefault}{\updefault}2}}}}}
\put(1312,1012){\makebox(0,0)[lb]{\smash{{{\SetFigFont{12}{14.4}{\rmdefault}{\mddefault}{\updefault}2}}}}}
\end{picture}
}
\end{center}
%%\vskip -60pt
\mycap{10}{mod $2$ representation of $\pi^{*} = 10 + 8 + 7 + 1$}    
\end{figure}

Note that the ordinary Ferrers graph representations do not possess this
invariance property. For example, if we conjugate the mod $1$ graph
in Fig.\ 9 we get the partition $5 + 5 + 4 + 4 + 4 + 3 + 1$,
which has repeated odd part $5$. 

Next, we define the $M_2$-rank of a partition as the
{\it largest row minus the number of rows of its mod $2$ graph}. 
It is easy to check that the $M_2$-rank of the partition,
$7 + 6 + 6 + 5 + 2$, depicted in Fig.\ 9, is equal to $4 - 5 = -1$,
while its rank is $7 - 5 = 2$. Also, it is straightforward to verify that
under conjugation the $M_2$-rank changes its sign, as does the
ordinary rank.

Let us define $\Etwid_r(q)$ as
\begin{equation}
\Etwid_r(q) = \sum_{n\ge1} \etwid_r(q) q^n,
\tag{5.2}
\end{equation}
where $\etwid_r(n)$ denotes the number of partitions of $n$
into distinct odd parts and unrestricted even parts such that
the $M_2$-rank $\ge r$. We assume that $\etwid_r(0)=0$.
We now show that
\begin{equation}
\Etwid_r(q) + \Etwid_{1-r}(q) + 1 = E(q),
\tag{5.3}
\end{equation}
and
\begin{equation}
\Etwid_r(q) = q^{2r+1} \left( \Etwid_{-1-r}(q) + 1\right), \quad r\ge0.
\tag{5.4}
\end{equation}
To prove \tagref{5.3} we will follow a well-trodden path and observe
that any nonempty partition counted by $E(q)$ whose $M_2$-rank
$\ge r$ is also counted by $\Etwid_r(q)$.
Any nonempty partition counted by $E(q)$ whose $M_2$-rank
$< r$ gives rise to a partition with $M_2$-rank
$\ge -1-r$, after conjugation. Thus, this conjugated partition
is counted by $\Etwid_{-1-r}(q)$ in \tagref{5.3}. Finally, the empty
partition is counted by 1 and $E(q)$ on the left and right sides
of \tagref{5.3}, respectively.

The proof of \tagref{5.4} requires modification of Dyson's transformation,
which we now proceed to describe. Let $\pi$ denote the mod $2$ 
graph of some partition counted by $\Etwid_r(q)$ in \tagref{5.4}.
Let $r + \ell(\pi)$ denote the length of the largest row of $\pi$,
and $h(\pi)$ denote the number of rows of $\pi$.
Clearly,
$$
h(\pi) \le \ell(\pi)
$$
for $\pi$ to have $M_2$-rank $\ge r$. Next, we remove the largest row
from $\pi$ to get a mod $2$ graph $\pitwid$. Conjugating $\pitwid$
we obtain $\pitwid^{*}$. Now, if the removed row represented the odd part
$2\ell + 2r - 1$, then we add to $\pitwid^{*}$ a new largest row
of length $\ell - 1$, representing the even part $2\ell - 2$.
On the other hand, if the removed row represented the even part $2\ell + 2r$,
then we add to $\pitwid^{*}$ a new largest row
of length $\ell$, representing the odd  part $2\ell - 1$.
These operations are illustrated in Fig. 11, 12 where the resulting partition
is denoted by $\pi'$.

\begin{figure}[h!]
\begin{center}
\setlength{\unitlength}{0.00066667in}
\begingroup\makeatletter\ifx\SetFigFont\undefined%
\gdef\SetFigFont#1#2#3#4#5{%
  \reset@font\fontsize{#1}{#2pt}%
  \fontfamily{#3}\fontseries{#4}\fontshape{#5}%
  \selectfont}%
\fi\endgroup%
{\renewcommand{\dashlinestretch}{30}
\begin{picture}(7437,1239)(0,-10)
\path(375,1212)(675,1212)(675,912)
	(375,912)(375,1212)
\path(675,1212)(975,1212)(975,912)
	(675,912)(675,1212)
\path(975,1212)(1275,1212)(1275,912)
	(975,912)(975,1212)
\put(475,1012){\makebox(0,0)[lb]{\smash{{{\SetFigFont{10}{12.0}{\rmdefault}{\mddefault}{\updefault}2}}}}}
\put(775,1012){\makebox(0,0)[lb]{\smash{{{\SetFigFont{10}{12.0}{\rmdefault}{\mddefault}{\updefault}2}}}}}
\put(1075,1012){\makebox(0,0)[lb]{\smash{{{\SetFigFont{10}{12.0}{\rmdefault}{\mddefault}{\updefault}2}}}}}
\path(1275,1212)(1575,1212)(1575,912)
	(1275,912)(1275,1212)
\path(375,612)(675,612)(675,312)
	(375,312)(375,612)
\path(975,912)(1275,912)(1275,612)
	(975,612)(975,912)
\path(675,912)(975,912)(975,612)
	(675,612)(675,912)
\path(375,912)(675,912)(675,612)
	(375,612)(375,912)
\thicklines
\path(1500,612)(2100,612)
\blacken\thinlines
\path(1980.000,522.000)(2100.000,612.000)(1980.000,702.000)(2016.000,612.000)(1980.000,522.000)
\path(2625,612)(2925,612)(2925,312)
	(2625,312)(2625,612)
\path(3225,912)(3525,912)(3525,612)
	(3225,612)(3225,912)
\path(2925,912)(3225,912)(3225,612)
	(2925,612)(2925,912)
\path(2625,912)(2925,912)(2925,612)
	(2625,612)(2625,912)
\thicklines
\path(3675,612)(4275,612)
\blacken\thinlines
\path(4155.000,522.000)(4275.000,612.000)(4155.000,702.000)(4191.000,612.000)(4155.000,522.000)
\path(5025,912)(5325,912)(5325,612)
	(5025,612)(5025,912)
\path(4725,912)(5025,912)(5025,612)
	(4725,612)(4725,912)
\path(4725,312)(5025,312)(5025,12)
	(4725,12)(4725,312)
\path(4725,612)(5025,612)(5025,312)
	(4725,312)(4725,612)
\path(6825,912)(7125,912)(7125,612)
	(6825,612)(6825,912)
\path(6525,912)(6825,912)(6825,612)
	(6525,612)(6525,912)
\path(6525,312)(6825,312)(6825,12)
	(6525,12)(6525,312)
\path(6525,612)(6825,612)(6825,312)
	(6525,312)(6525,612)
\path(6525,1212)(6825,1212)(6825,912)
	(6525,912)(6525,1212)
\path(6825,1212)(7125,1212)(7125,912)
	(6825,912)(6825,1212)
\path(7125,1212)(7425,1212)(7425,912)
	(7125,912)(7125,1212)
\thicklines
\path(5475,612)(6075,612)
\blacken\thinlines
\path(5955.000,522.000)(6075.000,612.000)(5955.000,702.000)(5991.000,612.000)(5955.000,522.000)
\put(1375,1012){\makebox(0,0)[lb]{\smash{{{\SetFigFont{10}{12.0}{\rmdefault}{\mddefault}{\updefault}1}}}}}
\put(0,762){\makebox(0,0)[lb]{\smash{$\pi$:}}}
\put(475,412){\makebox(0,0)[lb]{\smash{{{\SetFigFont{10}{12.0}{\rmdefault}{\mddefault}{\updefault}2}}}}}
\put(1075,712){\makebox(0,0)[lb]{\smash{{{\SetFigFont{10}{12.0}{\rmdefault}{\mddefault}{\updefault}1}}}}}
\put(775,712){\makebox(0,0)[lb]{\smash{{{\SetFigFont{10}{12.0}{\rmdefault}{\mddefault}{\updefault}2}}}}}
\put(475,712){\makebox(0,0)[lb]{\smash{{{\SetFigFont{10}{12.0}{\rmdefault}{\mddefault}{\updefault}2}}}}}
\put(2325,762){\makebox(0,0)[lb]{\smash{$\pitwid$:}}}
\put(2725,412){\makebox(0,0)[lb]{\smash{{{\SetFigFont{10}{12.0}{\rmdefault}{\mddefault}{\updefault}2}}}}}
\put(3325,712){\makebox(0,0)[lb]{\smash{{{\SetFigFont{10}{12.0}{\rmdefault}{\mddefault}{\updefault}1}}}}}
\put(3025,712){\makebox(0,0)[lb]{\smash{{{\SetFigFont{10}{12.0}{\rmdefault}{\mddefault}{\updefault}2}}}}}
\put(2725,712){\makebox(0,0)[lb]{\smash{{{\SetFigFont{10}{12.0}{\rmdefault}{\mddefault}{\updefault}2}}}}}
\put(4350,762){\makebox(0,0)[lb]{\smash{$\pitwid^{*}$:}}}
\put(5125,712){\makebox(0,0)[lb]{\smash{{{\SetFigFont{10}{12.0}{\rmdefault}{\mddefault}{\updefault}2}}}}}
\put(4825,712){\makebox(0,0)[lb]{\smash{{{\SetFigFont{10}{12.0}{\rmdefault}{\mddefault}{\updefault}2}}}}}
\put(4825,112){\makebox(0,0)[lb]{\smash{{{\SetFigFont{10}{12.0}{\rmdefault}{\mddefault}{\updefault}1}}}}}
\put(4825,412){\makebox(0,0)[lb]{\smash{{{\SetFigFont{10}{12.0}{\rmdefault}{\mddefault}{\updefault}2}}}}}
\put(6925,712){\makebox(0,0)[lb]{\smash{{{\SetFigFont{10}{12.0}{\rmdefault}{\mddefault}{\updefault}2}}}}}
\put(6625,712){\makebox(0,0)[lb]{\smash{{{\SetFigFont{10}{12.0}{\rmdefault}{\mddefault}{\updefault}2}}}}}
\put(6625,112){\makebox(0,0)[lb]{\smash{{{\SetFigFont{10}{12.0}{\rmdefault}{\mddefault}{\updefault}1}}}}}
\put(6625,412){\makebox(0,0)[lb]{\smash{{{\SetFigFont{10}{12.0}{\rmdefault}{\mddefault}{\updefault}2}}}}}
\put(6625,1012){\makebox(0,0)[lb]{\smash{{{\SetFigFont{10}{12.0}{\rmdefault}{\mddefault}{\updefault}2}}}}}
\put(6925,1012){\makebox(0,0)[lb]{\smash{{{\SetFigFont{10}{12.0}{\rmdefault}{\mddefault}{\updefault}2}}}}}
\put(7225,1012){\makebox(0,0)[lb]{\smash{{{\SetFigFont{10}{12.0}{\rmdefault}{\mddefault}{\updefault}2}}}}}
\put(6150,762){\makebox(0,0)[lb]{\smash{$\pi'$:}}}
\end{picture}
}
\end{center}
%%\vskip -60pt
\textsc{Fig 11}.\ \parbox[t]{4in}{Modification of Dyson's adjoint for
$\pi=7 + 5 + 2$ with $M_2$-rank $= 1 > r = 0$.}
\end{figure}

\begin{figure}[h!]
\begin{center}
\setlength{\unitlength}{0.00066667in}
\begingroup\makeatletter\ifx\SetFigFont\undefined%
\gdef\SetFigFont#1#2#3#4#5{%
  \reset@font\fontsize{#1}{#2pt}%
  \fontfamily{#3}\fontseries{#4}\fontshape{#5}%
  \selectfont}%
\fi\endgroup%
{\renewcommand{\dashlinestretch}{30}
\begin{picture}(7437,1239)(0,-10)
\path(375,1212)(675,1212)(675,912)
	(375,912)(375,1212)
\path(675,1212)(975,1212)(975,912)
	(675,912)(675,1212)
\path(975,1212)(1275,1212)(1275,912)
	(975,912)(975,1212)
\put(475,1012){\makebox(0,0)[lb]{\smash{{{\SetFigFont{10}{12.0}{\rmdefault}{\mddefault}{\updefault}2}}}}}
\put(775,1012){\makebox(0,0)[lb]{\smash{{{\SetFigFont{10}{12.0}{\rmdefault}{\mddefault}{\updefault}2}}}}}
\put(1075,1012){\makebox(0,0)[lb]{\smash{{{\SetFigFont{10}{12.0}{\rmdefault}{\mddefault}{\updefault}2}}}}}
\path(1275,1212)(1575,1212)(1575,912)
	(1275,912)(1275,1212)
\path(375,612)(675,612)(675,312)
	(375,312)(375,612)
\path(975,912)(1275,912)(1275,612)
	(975,612)(975,912)
\path(675,912)(975,912)(975,612)
	(675,612)(675,912)
\path(375,912)(675,912)(675,612)
	(375,612)(375,912)
\thicklines
\path(1500,612)(2100,612)
\blacken\thinlines
\path(1980.000,522.000)(2100.000,612.000)(1980.000,702.000)(2016.000,612.000)(1980.000,522.000)
\path(2625,612)(2925,612)(2925,312)
	(2625,312)(2625,612)
\path(3225,912)(3525,912)(3525,612)
	(3225,612)(3225,912)
\path(2925,912)(3225,912)(3225,612)
	(2925,612)(2925,912)
\path(2625,912)(2925,912)(2925,612)
	(2625,612)(2625,912)
\thicklines
\path(3675,612)(4275,612)
\blacken\thinlines
\path(4155.000,522.000)(4275.000,612.000)(4155.000,702.000)(4191.000,612.000)(4155.000,522.000)
\path(5025,912)(5325,912)(5325,612)
	(5025,612)(5025,912)
\path(4725,912)(5025,912)(5025,612)
	(4725,612)(4725,912)
\path(4725,312)(5025,312)(5025,12)
	(4725,12)(4725,312)
\path(4725,612)(5025,612)(5025,312)
	(4725,312)(4725,612)
\path(6825,912)(7125,912)(7125,612)
	(6825,612)(6825,912)
\path(6525,912)(6825,912)(6825,612)
	(6525,612)(6525,912)
\path(6525,312)(6825,312)(6825,12)
	(6525,12)(6525,312)
\path(6525,612)(6825,612)(6825,312)
	(6525,312)(6525,612)
\path(6525,1212)(6825,1212)(6825,912)
	(6525,912)(6525,1212)
\path(6825,1212)(7125,1212)(7125,912)
	(6825,912)(6825,1212)
\path(7125,1212)(7425,1212)(7425,912)
	(7125,912)(7125,1212)
\thicklines
\path(5475,612)(6075,612)
\blacken\thinlines
\path(5955.000,522.000)(6075.000,612.000)(5955.000,702.000)(5991.000,612.000)(5955.000,522.000)
\put(1375,1012){\makebox(0,0)[lb]{\smash{{{\SetFigFont{10}{12.0}{\rmdefault}{\mddefault}{\updefault}2}}}}}
\put(0,762){\makebox(0,0)[lb]{\smash{$\pi$:}}}
\put(475,412){\makebox(0,0)[lb]{\smash{{{\SetFigFont{10}{12.0}{\rmdefault}{\mddefault}{\updefault}2}}}}}
%%%
\put(1075,712){\makebox(0,0)[lb]{\smash{{{\SetFigFont{10}{12.0}{\rmdefault}{\mddefault}{\updefault}1}}}}}
\put(775,712){\makebox(0,0)[lb]{\smash{{{\SetFigFont{10}{12.0}{\rmdefault}{\mddefault}{\updefault}2}}}}}
\put(475,712){\makebox(0,0)[lb]{\smash{{{\SetFigFont{10}{12.0}{\rmdefault}{\mddefault}{\updefault}2}}}}}
\put(2325,762){\makebox(0,0)[lb]{\smash{$\pitwid$:}}}
\put(2725,412){\makebox(0,0)[lb]{\smash{{{\SetFigFont{10}{12.0}{\rmdefault}{\mddefault}{\updefault}2}}}}}
\put(3325,712){\makebox(0,0)[lb]{\smash{{{\SetFigFont{10}{12.0}{\rmdefault}{\mddefault}{\updefault}1}}}}}
\put(3025,712){\makebox(0,0)[lb]{\smash{{{\SetFigFont{10}{12.0}{\rmdefault}{\mddefault}{\updefault}2}}}}}
\put(2725,712){\makebox(0,0)[lb]{\smash{{{\SetFigFont{10}{12.0}{\rmdefault}{\mddefault}{\updefault}2}}}}}
\put(4350,762){\makebox(0,0)[lb]{\smash{$\pitwid^{*}$:}}}
\put(5125,712){\makebox(0,0)[lb]{\smash{{{\SetFigFont{10}{12.0}{\rmdefault}{\mddefault}{\updefault}2}}}}}
\put(4825,712){\makebox(0,0)[lb]{\smash{{{\SetFigFont{10}{12.0}{\rmdefault}{\mddefault}{\updefault}2}}}}}
\put(4825,112){\makebox(0,0)[lb]{\smash{{{\SetFigFont{10}{12.0}{\rmdefault}{\mddefault}{\updefault}1}}}}}
\put(4825,412){\makebox(0,0)[lb]{\smash{{{\SetFigFont{10}{12.0}{\rmdefault}{\mddefault}{\updefault}2}}}}}
\put(6925,712){\makebox(0,0)[lb]{\smash{{{\SetFigFont{10}{12.0}{\rmdefault}{\mddefault}{\updefault}2}}}}}
\put(6625,712){\makebox(0,0)[lb]{\smash{{{\SetFigFont{10}{12.0}{\rmdefault}{\mddefault}{\updefault}2}}}}}
\put(6625,112){\makebox(0,0)[lb]{\smash{{{\SetFigFont{10}{12.0}{\rmdefault}{\mddefault}{\updefault}1}}}}}
\put(6625,412){\makebox(0,0)[lb]{\smash{{{\SetFigFont{10}{12.0}{\rmdefault}{\mddefault}{\updefault}2}}}}}
\put(6625,1012){\makebox(0,0)[lb]{\smash{{{\SetFigFont{10}{12.0}{\rmdefault}{\mddefault}{\updefault}2}}}}}
\put(6925,1012){\makebox(0,0)[lb]{\smash{{{\SetFigFont{10}{12.0}{\rmdefault}{\mddefault}{\updefault}2}}}}}
\put(7225,1012){\makebox(0,0)[lb]{\smash{{{\SetFigFont{10}{12.0}{\rmdefault}{\mddefault}{\updefault}2}}}}}
\put(7525,1012){\makebox(0,0)[lb]{\smash{{{\SetFigFont{10}{12.0}{\rmdefault}{\mddefault}{\updefault}1}}}}}
\path(7425,1212)(7725,1212)(7725,912)
	(7425,912)(7425,1212)
\put(6150,762){\makebox(0,0)[lb]{\smash{$\pi'$:}}}
\end{picture}
}
\end{center}
%%\vskip -60pt
\textsc{Fig 12}.\ \parbox[t]{4in}{Modification of Dyson's adjoint for
$\pi=8 + 5 + 2$ with $M_2$-rank $= 1 > r = 0$.}
\end{figure}

Remarkably, regardless of whether the largest part of $\pi$ is even of odd
we have
\begin{equation}
|\pi'| = |\pi| - 2r - 1,
\tag{5.5}
\end{equation}
and
\begin{equation}
\mbox{$M_2$-rank}(\pi') \ge -1 - r.
\tag{5.6}
\end{equation}
It is easy to check that the map 
$\pi\to\pi'$ is reversible, except when $|\pi|=2r+1$. In the last case $\pi'$ is
empty. This concludes the proof of \tagref{5.4}.

Combining \tagref{5.3}, \tagref{5.4} we obtain
\begin{equation}
\Etwid_r(q) + q^{2r+1} \Etwid_{2+r}(q) = q^{2r+1} 
				\frac{ (-q;q^2)_\infty }{ (q^2;q^2)_\infty },
\quad r\ge 0.
\tag{5.7}
\end{equation}
Iteration of \tagref{5.7} yields
\begin{equation}
\Etwid_r(q) = \frac{ (-q;q^2)_\infty }{ (q^2;q^2)_\infty }
	 	\sum_{j\ge1} (-1)^{j-1} q^{2rj + j(2j-1)}, \quad r\ge 0.
\tag{5.8}
\end{equation}
Now \tagref{5.3} with $r=0$ states that
\begin{equation}
\Etwid_0(q) + \Etwid_1(q) + 1 = \frac{ (-q;q^2)_\infty }{ (q^2;q^2)_\infty }.
\tag{5.9}
\end{equation}
Thanks to \tagref{5.8} we may cast \tagref{5.9} in the form
\begin{equation}
1 = \frac{ (-q;q^2)_\infty }{ (q^2;q^2)_\infty }
	\sum_{j=-\infty}^{\infty} (-1)^j q^{j(2j+1)},
\tag{5.10}
\end{equation}
and so
\begin{equation}
\frac{ (q^2;q^2)_\infty }{ (-q;q^2)_\infty }
=        \sum_{j=-\infty}^{\infty} (-1)^j q^{j(2j+1)}.
\tag{5.11}
\end{equation}
Finally, replacing $q$ by $-q$ in \tagref{5.11} we obtain
the Gauss identity
\begin{equation}
\frac{ (q^2;q^2)_\infty }{ (q;q^2)_\infty }
=        \sum_{j=-\infty}^{\infty} q^{j(2j+1)} = \sum_{j\ge0} q^{T_j}.
\tag{5.12}
\end{equation}

Formula \tagref{5.8} implies that
\begin{align}
\Ehat_r(q) &= E(q) - \Etwid_{r+1}(q)\tag{5.13}\\
&= \frac{ (-q;q^2)_\infty}{ (q^2;q^2)_\infty }
   \sum_{j\ge0} (-1)^j q^{2rj + j(2j+1)},\quad r\ge0,
\nonumber
\end{align}
where $\Ehat_r(q)$ denotes the generating function for partitions into
distinct odd, and unrestricted even parts with $\mbox{$M_2$-rank}\le r$.
We now develop very different representations for $\Ehat_r(q)$.
To this end we decompose partitions counted by $\Ehat_r(q)$ into
even and odd parts. Let's assume that this decomposition gives
$\pi_1$ with $j$ distinct odd parts and $\pi_2$ with $i$ even parts.
Clearly, $\lambda(\pi_1) \le 2(i+j+r)-1$
and $\lambda(\pi_2) \le 2(i+j+r)$, and so, for $r\ge0$ we have
\begin{equation}
\Ehat_r(q) = \sum_{i,j\ge0}
q^{j + 2T_{j-1}} 
\qbinb{i+j+r}{j} q^{2i}
\qbinb{i+j+r-1+i}{i}.
\tag{5.14}
\end{equation}
Comparing \tagref{5.13} and \tagref{5.14}, we see that
\begin{equation}
\sum_{i,j\ge0} q^{j^2 + 2i}
\qbinb{i+j+r}{j} \qbinb{2i+j+r-1}{i}
=
\frac{ (-q;q^2)_\infty }{ (q^2;q^2)_\infty }
\sum_{j\ge0} (-1)^j q^{j(2j+1)}\left(q^{2r}\right)^j,\quad r\ge0.
\tag{5.15}
\end{equation}
Next, since 
\begin{equation}
\qbin{n+m}{n} = \frac{ (q^{1+m};q)_n }{ (q;q)_n },
\tag{5.16}
\end{equation}
we can rewrite \tagref{5.13} as
\begin{equation}
\sum_{i,j\ge0} q^{j^2 + 2i}
\frac{ (aq^{2i+2};q^2)_j (aq^{2i+2j};q^2)_i }
{ (q^2;q^2)_j (q^2;q^2)_i}
=
\frac{ (-q;q^2)_\infty }{ (q^2;q^2)_\infty }
\sum_{j\ge0} (-1)^j q^{2j^2+j} a^j,
\tag{5.17}
\end{equation}
where $a=q^{2r}$, $r\ge0$. Since the limit of the sequence
$\{q^{2r}\}$ is equal to zero, we may treat $a$ in \tagref{5.17}
as a free parameter.
%%%%%%%%%%%%%%%%%%%%%%%%%%%%%%%%%%%%%%%%%%%%%%%%%%%%%%%%%%%%%%%%%%%%%%%%%
%%INSERT #3
In Appendix B we discuss an alternative proof of 
\tagref{5.17}. This proof was communicated to us by George Andrews \cite{Y}.
%%%%%%%%%%%%%%%%%%%%%%%%%%%%%%%%%%%%%%%%%%%%%%%%%%%%%%%%%%%%%%%%%%%%%%%%%

In the past, fundamental as they are, modular representations have not
received the attention they deserve. Recently, Alladi \cite{AL}
used $2$-modular representations to provide an elegant
combinatorial bijection for a variant of G\"ollnitz's partition
theorem. However, in \cite{AL} partitions into only distinct odd parts
are considered, whereas here we allow even parts to appear with
possible repetition. 

In this regard, Alladi pointed out to us that Andrews \cite[ex.6, p.13]{A2}
used mod $2$ representations on the set of partitions treated here, subject
to the extra condition that no part $=1$, in order to establish a partition
theorem, which is equivalent to Cauchy's identity in the form:
\begin{equation}
\sum_{n\ge0} \frac{(-aq;q^2)_n}{(q^2;q^2)_n} t^n q^{2n}
= \frac{ (-atq^3;q^2)_\infty }{ (tq^2;q^2)_\infty }.
\tag{5.18}
\end{equation}
%%%%%%%%%%%%%%%%%%%%%%%%%%%%%%%%%%%%%%%%%%%%%%%%%%%%%%%%%%%%%%%%%%%%%%%%%
%%INSERT #4
Andrews's proof of the original Cauchy's identity with base $q$
(instead of base $q^2$ as above) may be found in \cite{Z}.
%%%%%%%%%%%%%%%%%%%%%%%%%%%%%%%%%%%%%%%%%%%%%%%%%%%%%%%%%%%%%%%%%%%%%%%%%
%%%%%END SECTION 5

\setcounter{equation}{0}
\section{Open questions}
\label{sec:6}

%%%%%%SECTION 6

In \cite{A3} Andrews proposed a dissection of a partition $\pi$
into successive Durfee squares with sizes
$n_1(\pi) \ge n_2(\pi) \ge n_3(\pi) \ge \cdots$. For example, the partition
$\pi$, depicted in Fig.\ 13, has two Durfee squares with sizes $n_1(\pi)=3$,
$n_2(\pi)=2$.

\begin{figure}[h!]
\begin{center}
\setlength{\unitlength}{0.00083333in}
\begingroup\makeatletter\ifx\SetFigFont\undefined%
\gdef\SetFigFont#1#2#3#4#5{%
  \reset@font\fontsize{#1}{#2pt}%
  \fontfamily{#3}\fontseries{#4}\fontshape{#5}%
  \selectfont}%
\fi\endgroup%
{\renewcommand{\dashlinestretch}{30}
\begin{picture}(1845,1860)(0,-10)
\path(33,312)(333,312)(333,12)
	(33,12)(33,312)
\path(333,612)(633,612)(633,312)
	(333,312)(333,612)
\path(33,612)(333,612)(333,312)
	(33,312)(33,612)
\path(333,912)(633,912)(633,612)
	(333,612)(333,912)
\path(33,912)(333,912)(333,612)
	(33,612)(33,912)
\path(333,1212)(633,1212)(633,912)
	(333,912)(333,1212)
\path(33,1212)(333,1212)(333,912)
	(33,912)(33,1212)
\path(933,1212)(1233,1212)(1233,912)
	(933,912)(933,1212)
\path(633,1212)(933,1212)(933,912)
	(633,912)(633,1212)
\path(333,1512)(633,1512)(633,1212)
	(333,1212)(333,1512)
\path(33,1512)(333,1512)(333,1212)
	(33,1212)(33,1512)
\path(933,1512)(1233,1512)(1233,1212)
	(933,1212)(933,1512)
\path(633,1512)(933,1512)(933,1212)
	(633,1212)(633,1512)
\path(1233,1512)(1533,1512)(1533,1212)
	(1233,1212)(1233,1512)
\path(1533,1812)(1833,1812)(1833,1512)
	(1533,1512)(1533,1812)
\path(1233,1812)(1533,1812)(1533,1512)
	(1233,1512)(1233,1812)
\path(933,1812)(1233,1812)(1233,1512)
	(933,1512)(933,1812)
\path(633,1812)(933,1812)(933,1512)
	(633,1512)(633,1812)
\path(333,1812)(633,1812)(633,1512)
	(333,1512)(333,1812)
\path(33,1812)(333,1812)(333,1512)
	(33,1512)(33,1812)
\thicklines
\path(33,1812)(933,1812)(933,912)
	(33,912)(33,1812)
\path(33,912)(633,912)(633,312)
	(33,312)(33,912)
\thinlines
\dashline{60.000}(33,1812)(1458,387)
\dashline{60.000}(33,912)(933,12)
\end{picture}
}
\end{center}
%%\vskip -60pt
\textsc{Fig 13}.\ \parbox[t]{5in}{Ferrers graph of $\pi=6 +5 + 4 + 2 + 2 + 1$.
This graph can be dissected into two Durfee squares of sizes $3$ and $2$.
$\mbox{$3$-rank}(\pi) = 2 - 1 = 1$.}
\end{figure}

Garvan \cite{Ga2} introduced a generalization of Dyson's rank
for partitions with at least $k-1$ successive Durfee squares.
He called this generalization the $k$-rank of a partition $\pi$.
The $k$-rank is defined as
\begin{equation}
\mbox{$k$-rank}(\pi) =
\matrix
&\mbox{the number of columns in the Ferrers graph of } \hfill\\
&\mbox{$\pi$ which lie to the right of the first Durfee square} \hfill\\
&\mbox{and whose length $\le n_{k-1}(\pi)$} \hfill\\
&\mbox{minus} \hfill\\
&\mbox{the number of parts of $\pi$ that lie below the} \hfill\\
&\mbox{$(k-1)$-th Durfee square.}\hfill
\endmatrix
\tag{6.1}
\end{equation}
For instance, the partition $\pi$ depicted in Fig.\ 13 has
$\mbox{$3$-rank}(\pi)=2-1=1$. Since any nonempty partition $\pi$
has at least one Durfee square we can easily infer that the $2$-rank
is the same as Dyson's rank.

Formula (1.10) in \cite{Ga2} implies that for $m\ge0$
\begin{equation}
FG_{k,m} (q) = \frac{1}{(q)_\infty} 
\sum_{j=1}^\infty (-1)^{j-1} q^{ \frac{j((2k-1)j-1)}{2} + mj},
\tag{6.2}
\end{equation}
where $FG_{k,m} (q)$ denotes the generating function for partitions
$\pi$ with at least $k-1$ successive Durfee squares and with
$\mbox{$k$-rank}(\pi) \ge m \ge 0$. Using \tagref{6.2} it is easy to
verify that
\begin{equation}
FG_{k,m}(q) + q^{k+m-1} FG_{k,2k-1+m}(q) = \frac{q^{k+m-1}}{ (q)_\infty }.
\tag{6.3}
\end{equation}
We note that \tagref{6.3} with $k=2$ becomes \tagref{1.24}.
Despite its speciously simple appearance the functional equation
\tagref{6.3} with $k > 2$ turned out to be very difficult to prove
in a combinatorial fashion. Perhaps the appropriate generalization
of Dyson's notion of rank-set may provide a key to a combinatorial
proof
of \tagref{6.3}.

We feel that it would be worthwhile to determine the precise
$q$-hypergeometric status of the new polynomial analogues of
Euler's pentagonal number theorem \tagref{1.30} and to explore
more general iteration schemes.
Finally, we would like to pose the problem of finding a natural bounded
extension of formulas \tagref{1.4}, \tagref{1.5}, and
\tagref{1.7}--\tagref{1.9}.
%%%%%END SECTION 6

\setcounter{equation}{0}
\section*{Appendix A}
%%\label{sec:6}
%%%%%%SECTION APPENDIX A
Here we give a direct proof of \tagref{3.17} which we restate as
\begin{equation}
\chat_{-k}(n) = \chat_k(n) + \delta_{n,1} \delta_{k,1}, \quad k\ge0,
\tag{A.1}
\end{equation}
with $\chat_{k}(n)$ denoting the number of partitions of $n$ with
crank $k$. It is easy to check that \tagref{A.1} holds for $n=0$, $1$.

To proceed further let us recall that $\lambda(\pi)$, $\mu(\pi)$
and $\nutwid(\pi)$ denote the largest part of a partition $\pi$,
the number of ones in $\pi$ and the number of parts of $\pi$ which are 
larger than $\mu(\pi)$, respectively. In addition, let $\gamma(\pi)$
be defined by
$$
\gamma(\pi) =
\begin{cases}
  p_1-p_2, &\mbox{if $\nutwid(\pi)\ne1$}, \\
  \lambda(\pi)-\mu(\pi)-1, &\mbox{if $\nutwid(\pi)=1$,}
\end{cases}
$$ 
where 
$\pi = p_1 + p_2 + p_3 + \cdots$ has
parts $p_1\ge p_2\ge p_3\ge\cdots$.
It is easy to check that \tagref{A.1} with $n>1$ is an immediate consequence 
of the following two propositions.

\proclaim{Proposition 1}
The number of partitions $\pi$ with $|\pi|=n>1$,
$\lambda(\pi)=\ell$,
and $\mu(\pi)=0$, equals the number of partitions $\pi'$ with
$|\pi'|=|\pi|$, $\nutwid(\pi')=0$, and $\mu(\pi')=\ell$.
\endproclaim

\proclaim{Proposition 2}
The number of partitions $\pi$ with $|\pi|=n>1$,
$\mu(\pi)=M>0$,
and $\nutwid(\pi)=N>0$,
equals the number of partitions $\pi'$ with
$|\pi'|=|\pi|$, 
$\mu(\pi)=N$,
and $\nutwid(\pi)=M$.
\endproclaim

To prove Proposition 1 we remove the largest row from the graph of $\pi$
and then add a vertical column representing $\lambda(\pi)$ ones to the
resulting graph. Let's call the resulting partition $\pi'$.
Obviously, $\mu(\pi')=\lambda(\pi)$, and
$\nutwid(\pi')=0$. 
Since $n>1$ the map $\pi\to\pi'$ is a bijection and the
result follows. \qed

The proof of Proposition 2 is more involved. Here we need to decompose
$\pi$ as indicated in Fig.\ 14 below.

%%%
%%% insert figure 14
\begin{figure}[h!]
\begin{center}
\setlength{\unitlength}{0.00083333in}
\begingroup\makeatletter\ifx\SetFigFont\undefined%
\gdef\SetFigFont#1#2#3#4#5{%
  \reset@font\fontsize{#1}{#2pt}%
  \fontfamily{#3}\fontseries{#4}\fontshape{#5}%
  \selectfont}%
\fi\endgroup%
{\renewcommand{\dashlinestretch}{30}
\begin{picture}(4287,3462)(0,-10)
\put(2175,2712){\ellipse{212}{212}}
\put(675,1662){\ellipse{212}{212}}
\path(1575,2112)(2475,2112)(2475,2262)
	(2625,2262)(2625,2412)(2925,2412)
	(2925,2712)(3375,2712)(3375,3012)
	(3375,3162)(1575,3162)
\texture{8101010 10000000 444444 44000000 11101 11000000 444444 44000000 
	101010 10000000 444444 44000000 10101 1000000 444444 44000000 
	101010 10000000 444444 44000000 11101 11000000 444444 44000000 
	101010 10000000 444444 44000000 10101 1000000 444444 44000000 }
\shade\path(3375,3162)(3525,3162)(3525,3012)
	(3375,3012)(3375,3162)
\path(3375,3162)(3525,3162)(3525,3012)
	(3375,3012)(3375,3162)
\shade\path(3525,3162)(3675,3162)(3675,3012)
	(3525,3012)(3525,3162)
\path(3525,3162)(3675,3162)(3675,3012)
	(3525,3012)(3525,3162)
\shade\path(3675,3162)(3825,3162)(3825,3012)
	(3675,3012)(3675,3162)
\path(3675,3162)(3825,3162)(3825,3012)
	(3675,3012)(3675,3162)
\shade\path(3825,3162)(3975,3162)(3975,3012)
	(3825,3012)(3825,3162)
\path(3825,3162)(3975,3162)(3975,3012)
	(3825,3012)(3825,3162)
\shade\path(3975,3162)(4125,3162)(4125,3012)
	(3975,3012)(3975,3162)
\path(3975,3162)(4125,3162)(4125,3012)
	(3975,3012)(3975,3162)
\shade\path(4125,3162)(4275,3162)(4275,3012)
	(4125,3012)(4125,3162)
\path(4125,3162)(4275,3162)(4275,3012)
	(4125,3012)(4125,3162)
\shade\path(1425,3162)(1575,3162)(1575,3012)
	(1425,3012)(1425,3162)
\path(1425,3162)(1575,3162)(1575,3012)
	(1425,3012)(1425,3162)
\shade\path(1425,3012)(1575,3012)(1575,2862)
	(1425,2862)(1425,3012)
\path(1425,3012)(1575,3012)(1575,2862)
	(1425,2862)(1425,3012)
\shade\path(1425,2862)(1575,2862)(1575,2712)
	(1425,2712)(1425,2862)
\path(1425,2862)(1575,2862)(1575,2712)
	(1425,2712)(1425,2862)
\shade\path(1425,2712)(1575,2712)(1575,2562)
	(1425,2562)(1425,2712)
\path(1425,2712)(1575,2712)(1575,2562)
	(1425,2562)(1425,2712)
\shade\path(1425,2562)(1575,2562)(1575,2412)
	(1425,2412)(1425,2562)
\path(1425,2562)(1575,2562)(1575,2412)
	(1425,2412)(1425,2562)
\shade\path(1425,2412)(1575,2412)(1575,2262)
	(1425,2262)(1425,2412)
\path(1425,2412)(1575,2412)(1575,2262)
	(1425,2262)(1425,2412)
\shade\path(1425,2262)(1575,2262)(1575,2112)
	(1425,2112)(1425,2262)
\path(1425,2262)(1575,2262)(1575,2112)
	(1425,2112)(1425,2262)
\shade\path(1425,2112)(1575,2112)(1575,1962)
	(1425,1962)(1425,2112)
\path(1425,2112)(1575,2112)(1575,1962)
	(1425,1962)(1425,2112)
\blacken\path(3495.000,3342.000)(3375.000,3312.000)(3495.000,3282.000)(3459.000,3312.000)(3495.000,3342.000)
\path(3375,3312)(3750,3312)
\path(3375,3312)(3750,3312)
\blacken\path(4155.000,3282.000)(4275.000,3312.000)(4155.000,3342.000)(4191.000,3312.000)(4155.000,3282.000)
\path(4275,3312)(3975,3312)
\path(4275,3312)(3975,3312)
\path(1425,1962)(1425,1662)(1125,1662)
	(975,1662)(975,1362)(675,1362)
	(675,1062)(375,1062)(375,1962)
\path(375,1962)(1425,1962)
\path(375,1962)(1425,1962)
\blacken\path(495.000,3192.000)(375.000,3162.000)(495.000,3132.000)(459.000,3162.000)(495.000,3192.000)
\texture{0 115111 51000000 444444 44000000 151515 15000000 444444 
	44000000 511151 11000000 444444 44000000 151515 15000000 444444 
	44000000 115111 51000000 444444 44000000 151515 15000000 444444 
	44000000 511151 11000000 444444 44000000 151515 15000000 444444 }
\path(375,3162)(1425,3162)
\path(375,3162)(1425,3162)
\blacken\path(1305.000,3132.000)(1425.000,3162.000)(1305.000,3192.000)(1341.000,3162.000)(1305.000,3132.000)
\blacken\path(405.000,3042.000)(375.000,3162.000)(345.000,3042.000)(375.000,3078.000)(405.000,3042.000)
\path(375,3162)(375,1962)
\blacken\path(345.000,2082.000)(375.000,1962.000)(405.000,2082.000)(375.000,2046.000)(345.000,2082.000)
\texture{8101010 10000000 444444 44000000 11101 11000000 444444 44000000 
	101010 10000000 444444 44000000 10101 1000000 444444 44000000 
	101010 10000000 444444 44000000 11101 11000000 444444 44000000 
	101010 10000000 444444 44000000 10101 1000000 444444 44000000 }
\shade\path(375,1062)(525,1062)(525,912)
	(375,912)(375,1062)
\path(375,1062)(525,1062)(525,912)
	(375,912)(375,1062)
\shade\path(375,912)(525,912)(525,762)
	(375,762)(375,912)
\path(375,912)(525,912)(525,762)
	(375,762)(375,912)
\shade\path(375,762)(525,762)(525,612)
	(375,612)(375,762)
\path(375,762)(525,762)(525,612)
	(375,612)(375,762)
\shade\path(375,612)(525,612)(525,462)
	(375,462)(375,612)
\path(375,612)(525,612)(525,462)
	(375,462)(375,612)
\shade\path(375,462)(525,462)(525,312)
	(375,312)(375,462)
\path(375,462)(525,462)(525,312)
	(375,312)(375,462)
\shade\path(375,312)(525,312)(525,162)
	(375,162)(375,312)
\path(375,312)(525,312)(525,162)
	(375,162)(375,312)
\shade\path(375,162)(525,162)(525,12)
	(375,12)(375,162)
\path(375,162)(525,162)(525,12)
	(375,12)(375,162)
\blacken\path(255.000,942.000)(225.000,1062.000)(195.000,942.000)(225.000,978.000)(255.000,942.000)
\path(225,1062)(225,12)
\path(225,1062)(225,12)
\blacken\path(195.000,132.000)(225.000,12.000)(255.000,132.000)(225.000,96.000)(195.000,132.000)
%%%%% math bits
\put(1575,1737){\makebox(0,0)[lb]{\smash{{{\SetFigFont{12}{14.4}{\rmdefault}{\mddefault}{\updefault}$\le M$}}}}}
\put(2135,2647){\makebox(0,0)[lb]{\smash{{{\SetFigFont{12}{14.4}{\rmdefault}{\mddefault}{\updefault}B}}}}}
\put(3825,3312){\makebox(0,0)[lb]{\smash{{{\SetFigFont{12}{14.4}{\rmdefault}{\mddefault}{\updefault}$\gamma$}}}}}
\put(75,2562){\makebox(0,0)[lb]{\smash{{{\SetFigFont{12}{14.4}{\rmdefault}{\mddefault}{\updefault}$N$}}}}}
\put(825,3312){\makebox(0,0)[lb]{\smash{{{\SetFigFont{12}{14.4}{\rmdefault}{\mddefault}{\updefault}$M$}}}}}
\put(0,462){\makebox(0,0)[lb]{\smash{{{\SetFigFont{12}{14.4}{\rmdefault}{\mddefault}{\updefault}$M$}}}}}
\put(635,1597){\makebox(0,0)[lb]{\smash{{{\SetFigFont{12}{14.4}{\rmdefault}{\mddefault}{\updefault}A}}}}}
\end{picture}
}
%%%%%%%%%%
\end{center}
\mycap{14}{\enskip Graph of $\pi$ in Proposition 2.}
\end{figure}
%%%
%%%

\noindent
Let us now remove from the graph of $\pi$ in Fig.\ 14 three pieces, namely, 
the vertical columns of height $M$, $N$ and the horizontal row of length 
$\gamma$. Next, we conjugate the resulting graph to get $\pitwid$. 
We now add three pieces to $\pitwid$ as indicated in Fig.\ 15 to get
$\pi'$. 

%%%
%%% insert figure 15
\begin{figure}[h!]
\begin{center}
\setlength{\unitlength}{0.00083333in}
\begingroup\makeatletter\ifx\SetFigFont\undefined%
\gdef\SetFigFont#1#2#3#4#5{%
  \reset@font\fontsize{#1}{#2pt}%
  \fontfamily{#3}\fontseries{#4}\fontshape{#5}%
  \selectfont}%
\fi\endgroup%
{\renewcommand{\dashlinestretch}{30}
\begin{picture}(3487,4335)(0,-10)
\put(875,2562){\ellipse{275}{275}}
\put(2000,3762){\ellipse{275}{275}}
\path(1673,3012)(1973,3012)(1973,3312)
	(1973,3462)(2274,3462)(2274,3762)
	(2574,3762)(2574,4062)(1674,4062)
\texture{8101010 10000000 444444 44000000 11101 11000000 444444 44000000 
	101010 10000000 444444 44000000 10101 1000000 444444 44000000 
	101010 10000000 444444 44000000 11101 11000000 444444 44000000 
	101010 10000000 444444 44000000 10101 1000000 444444 44000000 }
\shade\path(2725,4062)(2875,4062)(2875,3912)
	(2725,3912)(2725,4062)
\path(2725,4062)(2875,4062)(2875,3912)
	(2725,3912)(2725,4062)
\shade\path(2875,4062)(3025,4062)(3025,3912)
	(2875,3912)(2875,4062)
\path(2875,4062)(3025,4062)(3025,3912)
	(2875,3912)(2875,4062)
\shade\path(3025,4062)(3175,4062)(3175,3912)
	(3025,3912)(3025,4062)
\path(3025,4062)(3175,4062)(3175,3912)
	(3025,3912)(3025,4062)
\shade\path(3175,4062)(3325,4062)(3325,3912)
	(3175,3912)(3175,4062)
\path(3175,4062)(3325,4062)(3325,3912)
	(3175,3912)(3175,4062)
\shade\path(3325,4062)(3475,4062)(3475,3912)
	(3325,3912)(3325,4062)
\path(3325,4062)(3475,4062)(3475,3912)
	(3325,3912)(3325,4062)
\blacken\path(2695.000,4242.000)(2575.000,4212.000)(2695.000,4182.000)(2659.000,4212.000)(2695.000,4242.000)
\path(2575,4212)(2950,4212)
\path(2575,4212)(2950,4212)
\blacken\path(3355.000,4182.000)(3475.000,4212.000)(3355.000,4242.000)(3391.000,4212.000)(3355.000,4182.000)
\path(3475,4212)(3175,4212)
\path(3475,4212)(3175,4212)
\shade\path(2575,4062)(2725,4062)(2725,3912)
	(2575,3912)(2575,4062)
\path(2575,4062)(2725,4062)(2725,3912)
	(2575,3912)(2575,4062)
\shade\path(1500,4062)(1650,4062)(1650,3912)
	(1500,3912)(1500,4062)
\path(1500,4062)(1650,4062)(1650,3912)
	(1500,3912)(1500,4062)
\shade\path(1500,3912)(1650,3912)(1650,3762)
	(1500,3762)(1500,3912)
\path(1500,3912)(1650,3912)(1650,3762)
	(1500,3762)(1500,3912)
\shade\path(1500,3762)(1650,3762)(1650,3612)
	(1500,3612)(1500,3762)
\path(1500,3762)(1650,3762)(1650,3612)
	(1500,3612)(1500,3762)
\shade\path(1500,3612)(1650,3612)(1650,3462)
	(1500,3462)(1500,3612)
\path(1500,3612)(1650,3612)(1650,3462)
	(1500,3462)(1500,3612)
\shade\path(1500,3462)(1650,3462)(1650,3312)
	(1500,3312)(1500,3462)
\path(1500,3462)(1650,3462)(1650,3312)
	(1500,3312)(1500,3462)
\shade\path(1500,3312)(1650,3312)(1650,3162)
	(1500,3162)(1500,3312)
\path(1500,3312)(1650,3312)(1650,3162)
	(1500,3162)(1500,3312)
\shade\path(1500,3162)(1650,3162)(1650,3012)
	(1500,3012)(1500,3162)
\path(1500,3162)(1650,3162)(1650,3012)
	(1500,3012)(1500,3162)
\blacken\path(330.000,3942.000)(300.000,4062.000)(270.000,3942.000)(300.000,3978.000)(330.000,3942.000)
\path(300,4062)(300,3012)
\blacken\path(270.000,3132.000)(300.000,3012.000)(330.000,3132.000)(300.000,3096.000)(270.000,3132.000)
\blacken\path(420.000,4092.000)(300.000,4062.000)(420.000,4032.000)(384.000,4062.000)(420.000,4092.000)
\path(300,4062)(1500,4062)
\blacken\path(1380.000,4032.000)(1500.000,4062.000)(1380.000,4092.000)(1416.000,4062.000)(1380.000,4032.000)
\path(300,3012)(1500,3012)
\path(1350,3012)(1350,2113)(1200,2113)
	(1200,1962)(1050,1962)(1050,1662)
	(750,1662)(751,1212)(450,1212)
	(300,1212)(301,3012)
\blacken\path(180.000,1092.000)(150.000,1212.000)(120.000,1092.000)(150.000,1128.000)(180.000,1092.000)
\path(150,1212)(150,12)
\blacken\path(120.000,132.000)(150.000,12.000)(180.000,132.000)(150.000,96.000)(120.000,132.000)
\shade\path(300,1212)(450,1212)(450,1062)
	(300,1062)(300,1212)
\path(300,1212)(450,1212)(450,1062)
	(300,1062)(300,1212)
\shade\path(300,1062)(450,1062)(450,912)
	(300,912)(300,1062)
\path(300,1062)(450,1062)(450,912)
	(300,912)(300,1062)
\shade\path(300,912)(450,912)(450,762)
	(300,762)(300,912)
\path(300,912)(450,912)(450,762)
	(300,762)(300,912)
\shade\path(300,762)(450,762)(450,612)
	(300,612)(300,762)
\path(300,762)(450,762)(450,612)
	(300,612)(300,762)
\shade\path(300,612)(450,612)(450,462)
	(300,462)(300,612)
\path(300,612)(450,612)(450,462)
	(300,462)(300,612)
\shade\path(300,462)(450,462)(450,312)
	(300,312)(300,462)
\path(300,462)(450,462)(450,312)
	(300,312)(300,462)
\shade\path(300,312)(450,312)(450,162)
	(300,162)(300,312)
\path(300,312)(450,312)(450,162)
	(300,162)(300,312)
\shade\path(300,162)(450,162)(450,12)
	(300,12)(300,162)
\path(300,162)(450,162)(450,12)
	(300,12)(300,162)
\put(3025,4212){\makebox(0,0)[lb]{\smash{{{\SetFigFont{12}{14.4}{\rmdefault}{\mddefault}{\updefault}$\gamma$}}}}}
\put(1425,2862){\makebox(0,0)[lb]{\smash{{{\SetFigFont{12}{14.4}{\rmdefault}{\mddefault}{\updefault}$\le N$}}}}}
\put(825,4137){\makebox(0,0)[lb]{\smash{{{\SetFigFont{12}{14.4}{\rmdefault}{\mddefault}{\updefault}$N$}}}}}
\put(75,3462){\makebox(0,0)[lb]{\smash{{{\SetFigFont{12}{14.4}{\rmdefault}{\mddefault}{\updefault}$M$}}}}}
\put(-40,537){\makebox(0,0)[lb]{\smash{{{\SetFigFont{12}{14.4}{\rmdefault}{\mddefault}{\updefault}$N$}}}}}
\put(785,2497){\makebox(0,0)[lb]{\smash{{{\SetFigFont{12}{14.4}{\rmdefault}{\mddefault}{\updefault}B${}^*$}}}}}
\put(1910,3697){\makebox(0,0)[lb]{\smash{{{\SetFigFont{12}{14.4}{\rmdefault}{\mddefault}{\updefault}A${}^*$}}}}}
\end{picture}
}
\end{center}
\mycap{15}{\enskip Graph of $\pi'$ in Proposition 2.}
\end{figure}
%%%
%%%

\noindent
Clearly, $\mu(\pi')=N$, $\nutwid(\pi')=M$ and $|\pi'|=|\pi|$.
To finish the proof we observe that the map $\pi\to\pi'$ is
a bijection. \qed

Let us call the map employed in the proofs of Propositions 1 and 2
a pseudo-conjugation transformation. We say that a partition $\pi$
with $|\pi|>1$ is self-pseudo-conjugate if it remains invariant
under pseudo-conjugation. In addition, we say the
partitions $\pi=0$, $\pi=1$ are self-pseudo-conjugate.

It is well known that the number of self-conjugate partitions of $n$
equals the number of partitions into distinct odd parts. The generating
function for the last set of partitions is $(-q;q^2)_\infty$.
Remarkably, the same is true for self-pseudo-conjugate partitions,
as we now demonstrate. First, it is obvious that the partitions described in 
Proposition 1 are not self-pseudo-conjugate. Second, the partitions
$\pi$ in Proposition 2 are self-pseudo-conjugate only if $M=N$ and
the conjugate of sub-graph A in Fig.\ 14 is identical to
sub-graph B. Therefore, the generating function $\SPC(q)$ for
self-pseudo-conjugate partitions is
\begin{align}
\SPC(q) &= 1 + q + \sum_{\substack{M\ge1, \\ \gamma\ge0}}     
\frac{q^{M(M+1)+M+\gamma}}{(q^4;q^2)_{M-1}}\tag{A.2}\\
&=1 + q + \sum_{M\ge1}\frac{q^{M(M+1)+M}}{(1-q)(q^4;q^2)_{M-1}}
=(1+q) \sum_{M\ge0}\frac{q^{M(M+1)}}{(q^2;q^2)_{M}}q^M.
\nonumber
\end{align}
Making use of the Euler identity
\begin{equation}
\sum_{j\ge0} \frac{q^{j(j+1)}}{(q^2;q^2)_j} z^j = (-zq^2;q^2)_\infty,
\tag{A.3}
\end{equation}
we find that
\begin{equation}
\SPC(q) = (1+q) (-q^3;q^2)_\infty = (-q;q^2)_\infty,
\tag{A.4}
\end{equation}
as desired.
%%%%%%END SECTION APPENDIX A

\setcounter{equation}{0}
\section*{Appendix B}
%%\label{sec:6}
%%%%%%SECTION APPENDIX B
Here we describe an alternative proof of
\tagref{5.17} communicated to us by George Andrews \cite{Y}.
We begin by expanding the products
$(aq^{2i+2};q^2)_j$ and $(aq^{2i+2j};q^2)_i$
in \tagref{5.17} using the $q$-binomial theorem \cite[(3.3.6), p.36]{A2}
\begin{equation}
\sum_{n\ge0} q^{n^2-n} z^n \qbinb{L}{n}
= (-z;q^2)_L.
\tag{B.1}
\end{equation}
This way we obtain after preforming changes of summation variables the
following expression for the left side of \tagref{5.17}
\begin{equation}
\LHS\tagref{5.17}
=
\sum_{i,j,s,t\ge0}
(-a)^{s+t}
\frac{
q^{ (j+t)^2 + 2(i+s) + 2s(i+j+s+t) + s(s-1) + t^2+t + 2t(i+s) }
}
{
(q^2;q^2)_i
(q^2;q^2)_j
(q^2;q^2)_s
(q^2;q^2)_t
}.
\tag{B.2}
\end{equation}
Next, we use the Euler identity \tagref{A.3} along with another result
of Euler \cite[(2.2.5), p.19]{A2}
\begin{equation}
\sum_{n\ge0} \frac{z^n}{(q)_n} = \frac{1}{(z;q)_\infty}
\tag{B.3}
\end{equation}
to sum out the $j$ and $i$ variables in \tagref{B.2} to get

\begin{equation}
\LHS\tagref{5.17}
=
\sum_{s,t\ge0}
(-a)^{s+t}
\frac{q^{2t^2+4st+3s^2+s+t}}{(q^2;q^2)_s (q^2;q^2)_t}
\frac{(-q^{1+2s+2t};q^2)_\infty }{(q^{2+2s+2t};q^2)_\infty}.
\tag{B.4}
\end{equation}
Since

\begin{equation}
\frac{(-q^{1+2s+2t};q^2)_\infty }{(q^{2+2s+2t};q^2)_\infty}
= \frac{(-q;q^2)_\infty}{(q^2;q^2)_\infty}
\frac{(q^2;q^2)_{s+t}}{(-q;q^2)_{s+t}},
\tag{B.5}
\end{equation}
we can derive

\begin{equation}
\LHS\tagref{5.17}
=
 \frac{(-q;q^2)_\infty}{(q^2;q^2)_\infty}
\sum_{n\ge0} (-a)^n \frac{q^{2n^2+n}}{(-q;q^2)_n}
\sum_{s=0}^n q^{s^2} \qbinb{n}{s}.
\tag{B.6}
\end{equation}
Making use of \tagref{B.1} we can evaluate the inner sum in \tagref{B.6}
to get

\begin{equation}
\LHS\tagref{5.17}
=
\frac{(-q;q^2)_\infty}{(q^2;q^2)_\infty}
\sum_{n\ge0} (-a)^n q^{2n^2+n},
\tag{B.7}
\end{equation}
which is essentially \tagref{5.17}, as desired. 
%%%%%%END SECTION APPENDIX B

{\it Acknowledgements}.
We are grateful to Krishnaswami Alladi 
for interesting discussion and for his comments and suggestions,
and to George E. Andrews
for his help and encouragement.

{\it Note added}. In a recent paper, Warnaar \cite{W1} observed that
\tagref{1.30} with $n=1$ is a limiting case of
a rather non-trivial cubic summation formula
\begin{multline}
\sum_{k=0}^{\lfloor N/2 \rfloor} \frac{1-Aq^{4k}}{1-A}
\frac{(A,Aq^{N+1};q^3)_k}{(q,q^{-N};q)_k}
\frac{(q^{-N};q)_{2k}}{(Aq^{N+1};q)_{2k}}
\frac{(C,D;q)_k}{(Aq^3/C,Aq^3/D;q^3)_k}\, q^k \nonumber\\
=\begin{cases}\displaystyle
\frac{(Aq^3,q^{2-N}/C,q^{2-N}/D;q^3)_{\lfloor N/3 \rfloor}}
{(Aq^3/C,Aq^3/D,q^{2-N}/CD;q^3)_{\lfloor N/3 \rfloor}},
& N\not\equiv 2\pmod{3},\\[3mm]
0, & N\equiv 2\pmod{3},
\end{cases}
\tag{N.1}
\end{multline}
with $CD=Aq^{N+1}$.
More precisely, Warnaar replaced $A\to A^2$, $C\to CA$,
$D\to DA$ and let $A\to0$ in \tagref{N.1} to obtain
\begin{equation}
\sum_{k=0}^{\lfloor N/2 \rfloor} \frac{ (q^{-N};q)_{2k} }{ (q,q^{-N};q)_k } q^k =
\begin{cases}
(-1)^{\lfloor N/3 \rfloor} q^{-\frac{N(N-1)}{6}},& N\not\equiv 2\pmod{3}, \\[3mm]
0, &\mbox{otherwise},
\end{cases}
\tag{N.2}
\end{equation}
which is essentially \tagref{1.30} with $n=1$ and
$L=\lfloor (N+1)/3 \rfloor$. In \cite{W1} Warnaar established \tagref{N.1} by 
setting $p=0$ in his elliptic generalization of \tagref{N.1}
(Corollary 4.13 in \cite{W2}).
Here, we would like to point out that the cubic summation formula \tagref{N.1}
is a special case of the Gasper-Rahman transformation formula (3.19)
in \cite{GR2}. Indeed, setting $ac=d=A$ and $b=cq^{1+N}$ in this formula
we get
\begin{align}
&\sum_{k=0}^{\lfloor N/2 \rfloor} \frac{1-Aq^{4k}}{1-A}
\frac{(A,Aq^{N+1};q^3)_k}{(q,q^{-N};q)_k}
\frac{(q^{-N};q)_{2k}}{(Aq^{N+1};q)_{2k}}
\frac{(cq^{1+N},A/c;q)_k}{(Aq^{2-N}/c,cq^3;q^3)_k}\, q^k \nonumber\\
&\, = \frac{ (Aq;q)_N }{ (q^{-N};q)_N }
 \frac{(q^{1-2N};q^3)_N}{(Aq^{2-N};q^3)_N} \cdot \tag{N.3}\\
&\, \cdot {}_8W_7(Aq^{-1-N};Aq^{1+N},c,Ac^{-1}q^{-1-N},q^{1-N},q^{-N};q^3,q^3)
\nonumber
\end{align}
with ${}_{r+1}W_r(a_1;a_4,\dots,a_{r+1};q,z)$ defined as in \cite[(2.11.11)]{GR1}.
Note that for $N\equiv2\pmod{3}$, $N>0$
\begin{equation}
(q^{1-2N};q^3)_N = 0,
\tag{N.4}
\end{equation}
and, consequently, the right side of \tagref{N.3} becomes zero.
When $N\not\equiv2\pmod{3}$, the series ${}_8W_7$ in \tagref{N.3}
can be summed thanks to Jackson's $q$-Dougall's summation \cite[(II.22)]{GR1}.
As a result, we obtain \tagref{N.1} with $C=cq^{1+N}$ and $D=A/c$.

%%% HERE %%%%%%%%%%%%%%%%%%%%%%%%%%%%%%%%%%%%%%%%%%%%%%%%%%%%%%%%%%%%%%%%%%%

\end{document}